\documentclass[11pt,reqno]{amsart}
\usepackage{graphicx,color}
\usepackage{amssymb,amscd,latexsym,epsfig}
\usepackage[all]{xy}

\newcommand\PP{\mathbb P}
\newcommand\C{\mathbb C}
\newcommand\Q{\mathbb Q}
\newcommand\R{\mathbb R}
\newcommand\Z{\mathbb Z}
\newcommand\N{\mathbb N}

\newcommand\Hom{\operatorname{Hom}}
\newcommand\Aut{\operatorname{Aut}}

\newcommand{\ind}{\operatorname{Ind}}
\newcommand{\mult}{\operatorname{Mult}}
\newcommand{\out}{\operatorname{out}}
\newcommand{\trop}{\operatorname{trop}}

\newcommand{\bra}{\langle}
\newcommand{\ket}{\rangle}
\newcommand{\wed}{\wedge}

\newcommand{\m}{\textbf{m}}
\newcommand{\w}{\textbf{w}}

\makeatletter \@addtoreset{equation}{section} \makeatother

\newtheorem{thm}[equation]{Theorem}
\newtheorem{lem}[equation]{Lemma}
\newtheorem{cor}[equation]{Corollary}
\newenvironment{rmk}{\noindent\textbf{Remark}.}{\\}
\newenvironment{exa}{\noindent\textbf{Example}.}{\\}

\title[]{Universal covers and the GW/Kronecker correspondence}
\author{Jacopo Stoppa}
\date{}
\address{Department of Pure Mathematics and Mathematical Statistics, University of Cambridge, Wilberforce Road, Cambridge CB3 0WB, UK}
\email{J.Stoppa@dpmms.cam.ac.uk}
\begin{document}  
\begin{abstract} The tropical vertex is an incarnation of mirror symmetry found by Gross, Pandharipande and Siebert. It can be applied to $m$-Kronecker quivers $K(m)$ (together with a result of Reineke) to compute the Euler characteristics of the moduli spaces of their (framed) representations in terms of Gromov-Witten invariants (as shown by Gross and Pandharipande). In this paper we study a possible geometric picture behind this correspondence, in particular constructing rational tropical curves from subquivers of the universal covering quiver $\widetilde{K}(m)$. Additional motivation comes from the physical interpretation of $m$-Kronecker quivers in the context of quiver quantum mechanics (especially work of F. Denef).  
\end{abstract}
\maketitle
\section{Introduction}
The $m$-Kronecker quiver $K(m)$ is the bipartite quiver with $m$ edges directed from $v_1$ (the \emph{source}) to $v_2$ (the \emph{sink}):
\begin{center}
\centerline{
\xymatrix{{v_1} \ar@/^2pc/[rr]^{e_1} \ar@/^1pc/[rr]^{e_2} \ar@/_1pc/[rr]_{e_{m-1}} \ar@/_2pc/[rr]_{e_m}& \vdots & {v_2}}}
\end{center}
A stability condition (central charge) for its dimension vectors is specified by a pair of integers $(w_1, w_2)$. We will always refer to the choice $(w_1, w_2) = (1, 0)$. One can then form smooth, projective moduli spaces $\mathcal{M}^{s, B}_{K(m)}(d)$ for stable representations of $K(m)$ with dimension vector $d$ and a $1$-dimensional framing at $v_1$ (respectively $\mathcal{M}^{s, F}_{K(m)}(d)$ for a framing at $v_2$, see e.g. \cite{rein_smooth} for the general theory). By the results of Engel and Reineke \cite{rein_smooth} we have explicit formulae for the topological Euler characteristics $\chi(\mathcal{M}^{s, B}_{K(m)}(d))$ (and also for $F$-framings).

Here however we are interested in an alternative and rather surprising way of computing these Euler characteristics, using an incarnation of mirror symmetry known as the tropical vertex of Gross, Pandharipande and Siebert \cite{gps}. It turns out that computing the generating function
\[\sum_{k \geq 0} \chi(\mathcal{M}^{s, B}_{K(m)}(k a, k b)) x^{ka}y^{kb}\]
is equivalent to working out a Gromov-Witten theory for a family of algebraic surfaces.

Fix coprime positive integers $a, b$ and let $\PP(a, b, 1)$ be the weighted projective plane $(\C^3\setminus\{0\})/\C^*$, with action given by $\lambda\cdot(z_1, z_2, z_3) = (\lambda^a z_1, \lambda^b z_2, \lambda z_3)$. Its toric fan is given by the duals of the divisors $D_1, D_2, D_{\textrm{out}}$ cut out by $z_1, z_2, z_3$. We denote by $D^o_1, D^o_2, D^o_{\textrm{out}}$ the subschemes obtained by removing the three torus fixed points. Also choose length $m$ ordered partitions $P_a, P_b$ with sizes $|P_a| = k a, |P_b| = k b$ for an integer $k > 0$. Then the relevant invariants for us are 
\[N_{a, b}[(P_a, P_b)] \in \Q\]
counting rational curves in the weighted projective plane $\PP(a, b, 1)$ which pass through $m$-tuples of distinct points $x^1_1,\dots, x^1_m$ on $D^o_1$, respectively $x^2_1, \dots, x^2_m$ on $D^o_2$, with multiplicities specified by $P_a, P_b$ and which are tangent to $D^o_{\textrm{out}}$ to order $k$. As an example $N_{1,1}(2 + 1, 1+1+1) = 3$ counts plane rational cubics with a prescribed node which pass through $4$ other prescribed points, and with $D_{\textrm{out}}$ an inflectional tangent. We refer to \cite{gps} Sections 0.4 and 6.4  for precise definitions and further examples. The numbers $N_{a, b}[(P_a, P_b)]$ are well defined and independent of the choice of points.

The GW/$m$-Kronecker correspondence is the identity in $\mathbb{Q}[[x,y]]$
\[\exp\left(\sum_{k \geq 1}  \sum_{|P_a| = k a, |P_b| = k b} k N_{a, b}[(P_a, P_b)]  x^{k a} y^{k b}\right) \hskip4cm\]
\[\hskip2cm = \left(1 + \sum_{k \geq 1}\chi(\mathcal{M}^{s, B}_{K(m)}(k a, k b)) x ^{k a} y^{k b}\right)^{\frac{m}{a}}\]
\begin{equation}\label{KronGW}\hskip5cm = \left(1 + \sum_{k \geq 1}\chi(\mathcal{M}^{s, F}_{K(m)}(k a, k b)) x ^{k a} y^{k b}\right)^{\frac{m}{b}},\end{equation}
(summing over length $m$ ordered partitions $P_a, P_b$), first written down explicitly by Gross and Pandharipande \cite{gp} Corollary 3.

Gross and Pandharipande \cite{gp} Section 3.5 and Reineke \cite{rein} Section 6 have asked if there is actually a correspondence between curves and representations underlying the equality \eqref{KronGW}. In particular one could ask how to costruct a rational curve starting from a given framed representation of $K(m)$. 

This question was the original motivation for writing this paper. We hoped initially that it would be possible to construct a rational \emph{tropical} curve starting from a suitable framed representation of the \emph{universal covering} of the quiver, $\widetilde{K}(m)$ (due to Reineke and Weist). Our hope was motivated by the case of the standard Kronecker quiver $K(2)$, where we will see that this is roughly true. By the results of Weist \cite{weist} (see Theorem \ref{weistThm} and \eqref{weist} below) passing to $\widetilde{K}(2)$ is the same as localising with respect to the natural $(\C^*)^2$-action, so the Euler characteristics can be computed already on $\widetilde{K}(2)$. On the other hand, the GW invariants $N_{a, b}[(P_a, P_b)]$ do arise from certain tropical counts $N^{\trop}(\w)$, see \cite{gps} Theorem 3.4 and Proposition 5.3 (summarized in Theorem \ref{tropicalThm} below). 

When $m \geq 3$ however this approach becomes problematic and we are not able to construct a single rational tropical curve from a given representation. What we do instead is roughly the following. For a finite subquiver $Q \subset \widetilde{K}(m)$ and a ``perturbative" parameter $k \geq 1$ we construct a whole (finite) set of rational tropical curves $\mathcal{S}_{Q, k}$. By Weist's Theorem, representations of $Q$ embed in representations of $K(m)$, so we may think of the assignment $Q \mapsto \mathcal{S}_{Q, k}$ as a refinement of the construction in \cite{gps}, where a very similar set of curves arises simply from $K(m)$. Counting the curves in $\mathcal{S}_{Q,k}$ for $k \gg 1$ which satisfy some constraints, with a suitable weight (adapted from the usual tropical multiplicity), expresses the contribution of $Q$ to $\chi(\mathcal{M}^{s, B}_{K(m)}(d))$ for some dimension vector $d$. Constraints on the number of ``legs" of our tropical curves correspond to constraints on the dimension vector $d$. We denote these ad hoc counts $N^{\trop}_Q(\w)$. We do not claim that they are genuine tropical invariants, i.e. independent of a crucial choice made in their construction. But via the GW/Kronecker correspondence, we find a posteriori a way to think of $N^{\trop}_Q(\w')$ as the contribution of $Q$ to some genuine invariants $N^{\trop}(\w)$.   

Here is the plan of the paper. We collect the necessary preliminary notions and results in section \ref{prelim}. In section \ref{preCovers} we discuss universal covering quivers and Weist's Theorem. In section \ref{preRein} we introduce the cornerstone of our approach, Reineke's Theorem \ref{reinThm}. The genuine tropical counts $N^{\trop}(\w)$ and their connection to the GW invariants $N_{a, b}[(P_a, P_b)]$ are discussed in section \ref{preTropical}. We first present the construction of tropical curves from subquivers of $\widetilde{K}(m)$ under some very strong assumptions, essentially restricting us to $\widehat{K}(2)$. We do this both because we think that $\widehat{K}(2)$ is a good example and because one has stronger results in this case. The construction takes up sections \ref{naiveDiagrams}, \ref{naiveTrees} and the first part of \ref{curveConstruction}. The construction which holds for general $m$ hinges on the factorization/deformation technique of \cite{gps}, explained in section \ref{baseRings}, and is done in section \ref{curveConstruction}. The numbers $N^{\trop}_{Q, k}(\w)$ are defined at the end of that section. In section \ref{main} we present our results connecting the tropical curves obtained from $Q$ with its contribution to $\chi(\mathcal{M}^{s, B}_{K(m)}(d))$, see especially Corollary \ref{contribSimple} and Corollary \ref{contribGen}. Very important additional motivation for the present work came from the paper of F. Denef \cite{denef}, so we include in section \ref{denef} some remarks about quiver quantum mechanics. However our grasp of the necessary physical background is very limited, and our discussion will be hardly satisfactory to the experts. The reader would be well-advised to consult \cite{denef}.\\
\textbf{Acknowledgements.} This is an application of some of the ideas in \cite{gps}, \cite{rein} and \cite{weist}. It was motivated by conversations with So Okada and Thorsten Weist, and I take this opportunity to thank them. I am also grateful to Hiraku Nakajima, Markus Reineke, Elisa Tenni and Richard Thomas, as well as to RIMS, Kyoto and Trinity College, Cambridge.
\section{Preliminary notions and results}\label{prelim}
\subsection{Universal covering quivers}\label{preCovers} Let $Q$ be a quiver without closed loops, with vertices $Q_0$ and edges $Q_1$. The algebraic torus $T := (\C^*)^{|Q_1|}$ acts on the affine spaces of representations $\textrm{Rep}_Q(d)$ for $d \in \N Q_0$, by scaling the linear maps in a representation. Let us write $\textrm{X}(T) := \Hom(T, \C) \cong \Z Q_1$, the character group of $T$. The \emph{abelian universal covering quiver} of $Q$ (due to Reineke, see \cite{weist} Section 3) is the quiver $\widehat{Q}$ with vertices $\widehat{Q}_0 = Q_0 \times \textrm{X}(T)$ and arrows given by
\[(\alpha, \chi)\!: (i, \chi) \to (j, \chi + e_{\alpha})\]
for $\alpha\!: i \to j$ in $Q_1$ and $\chi \in \textrm{X}(T)$. Here $e_{\alpha}$ is the character corresponding to $\alpha \in Q_1$. We say that a dimension vector $\hat{d} \in \N \widehat{Q}_0$ is compatible with $d \in \N Q_0$ if $d_i = \sum_{\chi} \hat{d}_{i, \chi}$ for all $i \in Q_0$, and we write $\hat{d}\sim d$. There is an action of $\Z Q_1$ on $\widehat{Q}_0$ defined by $\lambda \cdot (i, \chi) = (i, \chi + \lambda)$, which extends to an action on dimension vectors $\N \widehat{Q}_0$ by linearity. In the following we will denote by $[\,\hat{d}\,]$ the equivalence class of $\hat{d} \in \N \widehat{Q}_0$. 

Suppose now that we fix a stability function $\Theta\!: \Z Q_0 \to \C$ and a dimension vector $d$ for which there are no strictly semistable objects (often we call such dimension vectors \emph{coprime}). Weist studied the fixed locus for the induced torus action on $\mathcal{M}^s_Q(d)$, proving the isomorphism 
\[(\mathcal{M}^s_Q(d))^{T} \cong \bigcup_{[\,\hat{d}\,]\sim d}\mathcal{M}^s_{\widehat{Q}}(\hat{d})\]
(see \cite{weist} Theorem 3.11). In turn each of the moduli spaces $\mathcal{M}^s_{\widehat{Q}}(\hat{d})$ admits a torus action, and this gives rise to a tower of fixed loci, described by representations of iterated abelian covering quivers. We may then ask if for a fixed $d$ this process stabilizes after a finite number of iterations, and what the iterated fixed locus looks like. Weist gave an answer in terms of the universal covering quiver of $Q$.

So let $W(Q)$ be the group of words on $Q$, generated by arrows and their formal inverses. The \emph{universal covering quiver} $\widetilde{Q}$ of $Q$ (see \cite{weist} Section 3.4) is the quiver with vertices $\widetilde{Q}_0 = Q_0 \times W(Q)$ and arrows given by 
\[(\alpha, w)\!: (i, w) \to (j, w\alpha)\]
for $\alpha\!: i \to j$ in $Q_1$ and $w \in W(Q)$. As in the abelian case we have the notion of a compatible dimension vector $\tilde{d} \in \N\widetilde{Q}_0$ for $d \in \N Q_0$, and an action of $W(Q)$ on $\N \widetilde{Q}_0$ given by $w'\cdot(i, w) = (i, w w')$, with equivalence classes $[\,\tilde{d}\,]$.
\begin{thm}[Weist \cite{weist} Theorem 3.16]\label{weistThm} For a fixed coprime dimension vector $d$ for $Q$ the iteration process stabilizes, and the iterated fixed locus can be identified with the disjoint union
\begin{equation}
\bigcup_{[\,\tilde{d}\,]\sim d} \mathcal{M}^s_Q(\tilde{d}).
\end{equation}
\end{thm}
In particular for topological Euler characteristics we get 
\begin{equation*}
\chi(\mathcal{M}^s_{Q}(d)) = \sum_{[\,\tilde{d}\,]\sim d} \chi(\mathcal{M}^{s}_{\widetilde{Q}}(\tilde{d})).
\end{equation*}

For our applications we need a small variant of this result, replacing the coprime condition on $d$ with the presence of a framing. We only state this for the Kronecker quivers $K(m)$. A \emph{B-framing} (respectively \emph{F-framing}) for a representation of $K(m)$ is the choice of a $1$-dimensional subspace $L \subset V_{v_1}$ (respectively $L \subset V_{v_2}$). Similarly $B$ or $F$ framings of a representation of $\widetilde{K}(m)$ are given by $1$-dimensional subspaces $L \subset V_{(v_1, w)}$ or $L \subset V_{(v_2, w)}$ for some $w \in W(K(m))$. There is a natural notion of stability for framed representations, which implies ordinary semistability (we refer to \cite{rein_smooth}). The framing rules out strictly semistable objects, so that we have smooth moduli spaces $\mathcal{M}^{s, B}_{K(m)}(d)$ and $\mathcal{M}^{s, (v_1, \alpha)}_{\widetilde{K}(m)}(\tilde{d})$ (and similarly for $F$ framings). Then one can check that the proof of Theorem \ref{weistThm} carries over to this framed case, giving for Euler characteristics
\begin{equation}\label{weist}
\chi(\mathcal{M}^{s, B}_{K(m)}(d)) = \sum_{[\,\tilde{d}\,]\sim d}\sum_{w}\chi(\mathcal{M}^{s,(v_1, w)}_{\widetilde{K}(m)}(\tilde{d})).
\end{equation}     

\begin{figure}[ht]
\centerline{\includegraphics[scale=.7]{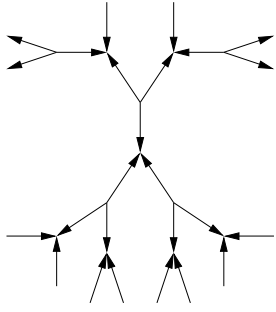}}
\caption{The universal covering quiver $\widetilde{K}(3)$.}
\label{unik3}
\end{figure}

We will often use the crucial fact (see \cite{weist} Remark 3.18) that the connected components of $\widetilde{K}(m)$ are given by infinite $m$-regular trees with an orientation.
\begin{figure}[ht]
\centerline{\includegraphics[scale=.7]{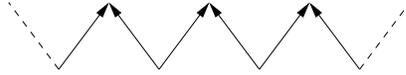}}
\caption{The universal abelian covering quiver $\widehat{K}(2)$.}
\label{k2}
\end{figure}
\begin{exa} The universal covering $\widetilde{K}(2)$ coincides with the universal abelian covering $\widehat{K}(2)$ (see Figure \ref{k2}).
\end{exa}
\begin{figure}[ht]
\centerline{\includegraphics[scale=.4]{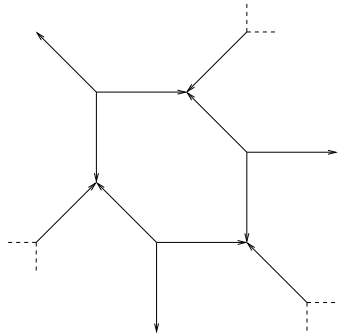}}
\caption{The universal abelian covering quiver $\widehat{K}(3)$.}
\label{k3}
\end{figure}
\begin{exa} The universal abelian covering quiver $\widehat{K}(3)$ is the infinite hexagonal quiver (see Figure \ref{k3}). The universal covering $\widetilde{K}(3)$ is obtained by opening up all the unoriented cycles in $\widehat{K}(3)$ (see Figure \ref{unik3}).
\end{exa}
\subsection{Reineke's theorem}\label{preRein} We start by fixing a finite subquiver $Q \subset \widetilde{K}(m)$. This is a bipartite quiver, i.e. every vertex is either a source or a sink. We label the sinks by $i_1, \dots, i_{s}$, the sources by $i_{s+1}, \dots, i_{s + S}$ (so there are $s$ sinks and $S$ sources). Notice that in particular $Q$ has no oriented (or indeed unoriented) cycles, so we can follow Reineke's convention and fix an order such that $i_k \to i_l \Rightarrow k > l$. For our purposes we also need that the order is \emph{minimal}, in the following sense: we label the sources mapping to $i_1$ by $i_{s + 1}, \ldots, i_{s + \ell_1}$, the sources mapping to $i_2$ by $i_{s + \ell_{1} + 1}, \ldots, i_{s + \ell_2}$, and so on.

A dimension vector $d$ has a reduction $\overline{d} \in \N K(m)_0 \cong \N \times \N$ given by 
\begin{equation}
\overline{d} = \big(\sum_{i > s} d_i, \sum_{i \leq s} d_i\big).
\end{equation}
We will write $\ind{d}, \ind(\overline{d})$ for the index of a dimension vector and its reduction, i.e. the unique positive integer $n$ such that $\frac{d}{n}$ (respectively $\frac{\overline{d}}{n}$) is primitive. Notice that we have $\ind(d) \leq \ind(\overline{d})$. The central charge $(1, 0) \in (\N K(m)_0)^{*}$ gives a notion of \emph{slope}, \begin{equation}
\mu(\overline{d}', \overline{d}'') = \frac{\overline{d}'}{\overline{d}' + \overline{d}''}.
\end{equation}
We fix a notion of slope for dimension vectors of $Q$ induced from the central charge $(1, 0)$ on $K(m)$, namely 
\begin{equation}\label{slope}
\mu(d) = \frac{\sum_{k > s} d_{k}}{\sum_k d_k}.
\end{equation} 
The set of dimension vectors with slope $\mu$ (together with the trivial representation) forms a subsemigroup $(\N Q_0)_{\mu} \subset \N Q_0$. 

The \emph{Euler form} is a bilinear form on $\Z Q_0$ defined by
\begin{equation}
e(d', d'') = \sum_{i \in Q_0} d'_i d''_i - \sum_{\alpha: i \to j} d'_i d''_j,
\end{equation} 
where the second sum is over all arrows from $i$ to $j$. We denote its skew-symmetrization by
\begin{equation}
\bra d', d''\ket = e(d', d'') - e(d'', d').
\end{equation}
\begin{rmk} The form $\bra\,\cdot, \cdot\ket$ is sometimes called the \emph{DSZ product} in physics terminology. Notice that in our case the product $\bra i_k, i_l\ket$ takes values in $\{0, \pm 1\}$. A possible source of confusion is that the skew-symmetrised Euler form is denoted by $\{\cdot, \cdot\}$ in Reineke's notation. 
\end{rmk}
A crucial role is played by a Poisson algebra modelled on $Q$, 
\[\mathcal{B} = (\C[[x_k]]_{k \in Q_0}, \{\,\cdot, \cdot\}),\]
with Poisson bracket generated by $\{ x_k, x_l \} = \bra k, l\ket x_k x_l$. For any dimension vector $d \in \N Q_0$ the Kontsevich-Soibelman Poisson automorphism $T_{d} \in \Aut(\mathcal{B})$ (a version of the operators appearing in \cite{ks} Section 1.4) is defined by
\begin{equation}
T_{d} (x_k) = x_k (1 + x^{d})^{\bra d, k \ket}.
\end{equation}
The fundamental object for us is the Poisson automorphism of $\mathcal{B}$ given by
\begin{equation}\label{product}
T_{i_1} \circ T_{i_2} \dots \circ T_{i_s} \circ T_{i_{s+1}} \circ \dots \circ T_{i_{s + S}}.
\end{equation}
By the general theory (see e.g. \cite{gps} Theorem 1.4) this can be written as a product of Poisson automorphisms attached to each rational nonegative slope, $\prod^{\leftarrow}_{\mu} \theta_{Q,\mu}$. The symbol $\leftarrow$ means we are writing factors in this product in the \emph{descending} slope order from left to right. Reineke showed that the Poisson automorphisms $\theta_{Q,\mu}$ can be computed in terms of the Euler characteristics of moduli spaces of stable framed representations of $Q$.
\begin{thm}[Reineke \cite{rein} Theorem 2.1]\label{reinThm} We have
\begin{equation}
\theta_{Q,\mu}(x_j) = x_j \prod_{i \in Q_0}(\theta_{Q, \mu,i}(x))^{\bra i, j\ket},
\end{equation}
where  
\begin{equation}
\theta_{Q, \mu,i}(x) = \sum_{d \in (\N Q_0)_{\mu}} \chi(\mathcal{M}^{s, i}_{Q}(d))\,x^d
\end{equation}
and $\mathcal{M}^{s, i}_{Q}(d)$ is the moduli space of stable representations of $Q$ (with respect to the choice of slope \eqref{slope}) with a $1$-dimensional framing at $i \in Q_0$.
\end{thm}
\begin{rmk} While we are only concerned with finite subquivers of $\widetilde{K}(m)$ we should make it clear that Reineke's result holds for general finite quivers without oriented cycles.
\end{rmk}
\subsection{Sorting diagrams}\label{naiveDiagrams} We will be concerned with an iterative process which sorts the factor of the fundamental product \eqref{product} in the opposite slope order, possibly introducing new factors at each step. This process is encoded by \emph{sorting diagrams}. Here we give a definition inspired by that of the \emph{scattering diagrams} appearing in \cite{gps} Section 1.4. For the sake of exposition, we initially make very strong assumptions about \eqref{product}, and give examples of sorting diagrams which these hold. The restrictive assumptions will be removed in the section \ref{baseRings} by working over more general base rings. 

So to a fixed product \eqref{product} we associate a unique sequence of \emph{sorting diagrams} $\mathfrak{S}^i, i \geq 0$. These are simply ordered sequences of group elements $\sigma^i_j \in \mathcal{B}$, 
\begin{equation}
\mathfrak{S}^i = (\sigma^i_1, \dots, \sigma^i_{\ell_i}).
\end{equation}
We set 
\begin{equation}
\mathfrak{S}^0 = (T_{i_1} , T_{i_2} , \dots T_{i_s} , T_{i_{s+1}} , \dots T_{i_{s + S}}).
\end{equation}
Notice that for all elements of $\mathfrak{S}^0$ we have a well defined notion of slope: since $\sigma^0_j = T_{i_j}$ we set $\mu(\sigma^0_j) = \mu(i_j)$. We define the $\mathfrak{S}^i$ for $i > 0$ inductively as follows. We move along the sequence $\mathfrak{S}^i$ starting from the left until we meet a pair of group elements with $\mu(\sigma^i_p) < \mu(\sigma^i_{p+1})$. We wish to commute $\sigma^i_p$ past the elements to its right until we meet again an element with smaller slope, $\mu(\sigma^i_{p}) \geq \mu(\sigma^i_{q})$. 

\textbf{Assumption 1.} Our first (very restrictive) assumption is that for $p + 1 \leq p' \leq q - 1$ we have 
\begin{equation}\label{assumption}
(\sigma^i_{p'})^{-1} \sigma^i_p \sigma^i_{p'} (\sigma^i_p)^{-1} = T_{d_{p'}} 
\end{equation} 
for some $d_{p'} \in \N Q_0$ (we follow the convention that $T_0 = 1$). Then we define $\mathfrak{S}^{i+1}$ by replacing the segment 
\[(\sigma^i_p, \dots, \sigma^i_{q-1})\] 
in $\mathfrak{S}^i$ by 
\[(\sigma^i_{p+1}, T_{d_{p+1}}, \sigma^i_{p+2}, T_{d_{p+2}}, \dots, T_{d_{q-1}}, \sigma^i_p).\] 
In particular $\mu(\sigma^{i+1}_j)$ is well defined for all $\sigma^{i+1}_j \in \mathfrak{S}^{i+1}$. 

\textbf{Assumption 2.} The sequence of sorting diagrams $\mathfrak{S}^{i}$ stabilizes for $i \gg 1$. We write $\mathfrak{S}^{\infty}$ for the stable sorting diagram.

\textbf{Assumption 3.} Operators in $\mathfrak{S}^{\infty}$ with the same slope commute.

If $\sigma, \tau$ are two operators in $\mathfrak{S}^{\infty}$ with $\sigma$ preceeding $\tau$, we write $\sigma \prec \tau$.

The following simple lemma is enough to effectively compute sorting diagrams, under their present (restrictive) definition.
\begin{lem}\label{simpleCommutator} If $\bra d, e\ket = 0$ then $ T_d \circ T_e = T_e \circ T_d$; and if $\bra d, e \ket = 1$ then $T_{d}\circ T_{e} = T_e \circ T_{d + e} \circ T_d$. 
\end{lem}
\begin{proof} Both equalities can be checked by direct computation, the second is the ``pentagon identity" \cite{ks} Section 1.4.
\end{proof}
\textbf{Example.} Consider a localization quiver $Q_1 \subset \widehat{K}(2)$ with underlying graph given by
\begin{center}\label{simplest}
\centerline{
\xymatrix{               &  {i_1}\\
     {i_3} \ar[ur]\ar[dr]&       \\
                         & {i_2}}}
\end{center}
The sorting diagrams $\mathfrak{S}^i$ stabilize for $i \geq 2$ and we find
\begin{align*}
\mathfrak{S}^0 &= (T_{i_1}, T_{i_2}, T_{i_3}),\\
\mathfrak{S}^1 &= (T_{i_1}, T_{i_3}, T_{i_2 + i_3}, T_{i_2}),\\
\mathfrak{S}^2 &= (T_{i_3}, T_{i_1 + i_3}, T_{i_2 + i_3}, T_{i_1 + i_2 + i_3}, T_{i_1}, T_{i_2}). 
\end{align*}
\textbf{Example.} Consider a localization quiver $Q_2 \subset \widehat{K}(2)$ with underlying graph given by
\begin{center}\label{simplest2}
\centerline{
\xymatrix{                
     {i_4} \ar[r]\ar[dr] &  {i_1}\\
                         &  {i_2}\\
     {i_5} \ar[ur]\ar[r] &  {i_3}}}                           
\end{center}
One can check that the sorting diagrams stabilize for $i \geq 6$. The first few are given by
\begin{align*}
\mathfrak{S}^0 &= (T_{i_1}, T_{i_2}, T_{i_3}, T_{i_4}, T_{i_5}),\\
\mathfrak{S}^1 &= (T_{i_1}, T_{i_2}, T_{i_4}, T_{i_5}, T_{i_3 + i_5}, T_{i_3}),\\
\mathfrak{S}^2 &= (T_{i_1}, T_{i_4}, T_{i_2 + i_4}, T_{i_5}, T_{i_2 + i_5},   T_{i_3 + i_5}, T_{i_2 + i_3 +i_5}, T_{i_2}, T_{i_3})\\
\mathfrak{S}^3 &= (T_{i_4}, T_{i_1 + i_4}, T_{i_2 + i_4}, T_{i_1 + i_2 + i_4},   T_{i_5}, T_{i_2 + i_5}, T_{i_3 + i_5}, T_{i_2 + i_3 + i_5}, T_{i_1}, T_{i_2},    T_{i_3})\\
\mathfrak{S}^4 &= (T_{i_4}, T_{i_1 + i_4}, T_{i_2 + i_4}, T_{i_5}, T_{i_1 + i_2 + i_4 + i_5}, T_{i_2 + i_5}, T_{i_3 + i_5}, T_{i_1 + i_2 + i_3 + i_4 + i_5},\\
& \,\,\,\,\,\,\,\,\,\,\,T_{i_1 + i_2 + i_4}, T_{i_2 + i_3 + i_5}, T_{i_1}, T_{i_2}, T_{i_3}).
\end{align*}
We also give an example where the present naive definition of sorting diagrams breaks down (i.e. the sorting diagram is undefined for some finite $i > 1$).\\
\textbf{Example.} Consider a localization quiver $Q_3 \subset \widetilde{K}(3)$ with underlying graph given by 
\begin{center}
\centerline{
\xymatrix{               
     {i_5} \ar[r]\ar[dr] &  {i_1} &\\                       
                         &  {i_2} & \ar[l] {i_7} \ar[d]\\
     {i_6} \ar[ur]\ar[r] &  {i_3} & {i_4}}}                         
\end{center}
A tedious but straighforward computation using Lemma \ref{simpleCommutator} shows that 
\begin{equation}
\mathfrak{S}^{6} = (\dots, T_{i_1 + i_2 + i_3 + i_5 + i_6}, T_{i_7}, T_{i_1 + i_2 + i_5 + i_7}, T_{i_2 + i_3 + i_6 + i_7}, T_{i_1 +2 i_2 + i_3 + i_5 + i_6 + i_7},\dots).
\end{equation}   
Set 
\begin{align*}
\xi  &= i_1 + i_2 + i_3 + i_5 + i_6,\\
\eta &= i_1 +2 i_2 + i_3 + i_5 + i_6 + i_7.
\end{align*}
The slopes of the elements in the displayed segment are $\{\frac{2}{5}, 1, \frac{1}{2}, \frac{1}{2}, \frac{3}{7}\}$. Therefore to compute $\mathfrak{S}^{7}$ we should commute $T_{\xi}$ past all the other elements in this segment. This works initially since
\begin{align*}
\bra \xi, i_7\ket &= 0,\\
\bra \xi, i_1 + i_2 + i_5 + i_7\ket &= 0,\\ 
\bra \xi, i_2 + i_3 + i_6 + i_7\ket &= 0,
\end{align*}
but at the last step we find
\begin{equation*}
\bra \xi, \eta \ket = -1.
\end{equation*}
We claim that the product $T^{-1}_{\eta} T_{\xi} T_{\eta} T^{-1}_{\xi}$ is not given by a single Poisson automorphism $T_d$ as in \eqref{assumption}. To see this consider the action of $T_{\xi}, T_{\eta}$ on the subalgebra generated by variables of the form $x^{a \xi + b \eta}$. We write $\overline{T}_{\xi}, \overline{T}_{\eta}$ for these restricted operators. Products of the form $\overline{T}^{\,-1}_{\eta} \overline{T}_{\xi} \overline{T}_{\eta} \overline{T}^{\,-1}_{\xi}$ are studied in \cite{ks} Section 1.4. It is shown there that there is a slope-ordered expansion
\begin{equation}
\overline{T}_{\xi} \overline{T}_{\eta} = \overline{T}_{\eta} \prod^{\to}_{a, b} \overline{T}^{\Omega(a, b)}_{a\eta + b\xi}\,\overline{T}_{\xi}
\end{equation} 
for certain $\Omega(a, b) \in \mathbb{Q}$, which are nonzero for \emph{infinitely many} values of $(a, b)$ as soon as $\bra \xi, \eta \ket \leq -1$ or $\bra \xi, \eta \ket \geq 2$. As observed by Kontsevich and Soibelman, a closed formula for the $\Omega(a, b)$ for $\bra \xi, \eta \ket \leq -1$ is not yet known. The first few terms are given by
\begin{equation}
\overline{T}_{\xi} \overline{T}_{\eta} \approx \overline{T}_{\eta} \overline{T}^{\,-1}_{3\xi + \eta} \overline{T}_{2\xi + \eta} \overline{T}^{2}_{3\xi + 2\eta} \overline{T}^{\,-1}_{\xi + \eta} \overline{T}^{\,-2}_{2\xi + 2\eta} \overline{T}^2_{2\xi + 3\eta} \overline{T}_{\xi + 2\eta} \overline{T}^{\,-1}_{\xi + 3 \eta} \overline{T}_{\xi}.
\end{equation}

This is enough to show that our present definition of sorting diagrams is too weak in general. There is however a special case when it (almost) works, that of the (abelian) universal covering quiver $\widetilde{K}(2)\cong \widehat{K}(2)$.
\begin{lem}\label{special1} Let $Q \subset \widehat{K}(2)$ be a finite subquiver. Then the sorting diagrams $\mathfrak{S}^{i}$ exist for all $i \geq 0$ and stabilize for $i \gg 1$ to a stable diagram $\mathfrak{S}^{\infty}$. Moreover operators in $\mathfrak{S}^{\infty}$ having the same slope $\mu \neq \frac{1}{2}$ commute.
\end{lem}
\begin{proof} Let $d$ be any dimension vector for $Q$. Recall that the moduli space of stable representations $\mathcal{M}^s_Q(d)$, when not empty, has dimension 
\begin{equation}\label{dim}
1 - e(d, d).
\end{equation} 
In our case of $Q \subset \widehat{K}(2)$ it is not hard to show, by induction, that the Euler form $e(d, d)$ is a positive definite quadratic form. It follows that the moduli space $\mathcal{M}^s_Q(d)$ must be empty for all but finitely many $d$. In fact we even have $e(d, d) > 1$ if $d_{i} > 1$ for some $i \in Q_0$.    

Suppose then that to form some sorting diagram $\mathfrak{S}^i$ for $Q$ we must commute some operator $T_{d'}$ past $T_{d''}$ with $\bra d', d''\ket \notin \{0, 1\}$. Then according to \cite{ks} Section 1.4, the ordered product expansion for $T^{-1}_{d''}\circ T_{d'} \circ T_{d''}\circ T^{-1}_{d'}$ must contain infinitely many factors $T_{d}$ with distinct slopes, $\mu(d') < \mu(d) < \mu(d'')$. Therefore by Reineke's Theorem \ref{reinThm} there must exist nonempty moduli spaces of ($B$ or $F$) framed representations for infinitely many $d$ with distinct slopes. Since framed stability implies semistability, the moduli spaces $\mathcal{M}^{ss}_Q(d)$ must by nonempty for infinitely many $d$ with distinct slopes. By the Jordan-Holder filtration, there also exist infinitely many distinct nonempty moduli spaces $\mathcal{M}^s_Q(d^*)$ for some dimension vectors $d^*$, which is a contradiction. This shows that the $\mathfrak{S}^i$ exist for all $i \geq 1$. Suppose that the $\mathfrak{S}^i$ do not stabilize for $i \gg 1$. Then for $i \gg 1$ the diagram $\mathfrak{S}^i$ must contain an operator $T_d$ with 
\begin{equation}
d = \sum_{i \in A}d_i + \sum_{j \in B}d_j
\end{equation} 
where $A, B$ are distinct subsets of $Q_0$ with $A \cap B \neq \emptyset$. In particular $d$ is coprime, so the appeareance of $T_d$ implies that $\mathcal{M}^{s}_Q(d)$ in nonempty. But for some $i \in Q_0$ we have $d_i = 2$, which again contradicts the dimension formula \eqref{dim}.

Notice that we have actually proved a much stronger result: all the operators appearing in $\mathfrak{S}^{\infty}$ must be of the form $T_d$ for $d = \sum_{i \in A} i$, where $A$ is some nonempty subset of $Q_0$. Moreover, by (framed) stability, the support of the dimension vector $d$ must be connected. 

It follows that for two dimension vectors $d', d''$ with the same slope $\mu \neq \frac{1}{2}$, the number of sources (respectively sinks) in $d'$ and $d''$ is the same (this is clearly \emph{not} true for $\mu = \frac{1}{2}$). Suppose $d' = \sum_{i \in A} i$, $d'' = \sum_{j \in B} j$. We can easily reduce to the case when $A \cap B = \emptyset$, and in this latter case there are obviously no arrows from the support of $d'$ to that of $d''$, so $\bra d', d''\ket = 0$. Therefore $T_{d'}, T_{d''}$ commute.  
\end{proof}
\begin{rmk} It is easy to show by example that we may have $[T_{d'}, T_{d''}] \neq 1$ when $\mu(d') = \mu(d'') = \frac{1}{2}$. This difficulty is related to the fact that the moduli spaces of stable framed representations of dimension vector $(n,n)$ are nonempty \emph{for all } $n$, namely $\mathcal{M}^{s, B}_{K(2)}(n, n) \cong \PP^{n}$. We do not address this problem in the special example of $K(2)$: it will be solved automatically when working over the more general base rings of section \ref{baseRings}. 

Let us also write down for later use the generating series of the Euler characteristics of stable $B$-framed representations for $K(2)$. In fact the only possible dimension vectors are those proportional to one of $(a, a+1), (1, 1)$ or $(a+1, a)$ for $a \geq 1$ (see e.g. \cite{gp} Lemma 2.3), and we have (see e.g. \cite{gp} Section 1.4 and Theorem 1)
\begin{align}
B_{a, a+1} &= \sum_{k\geq 1} \chi(\mathcal{M}^{s, B}_{K(2)}) x^{k a}y^{k(a+1)} = (1 + x^a y^{a + 1})^{a},\label{genSeries}\\
B_{1, 1} &= \sum_{k \geq 1}\chi(\mathcal{M}^{s, B}_{K(2)}) x^{k}y^{k} = (1 - x y)^{-2},\\
B_{a+1, a} &= \sum_{k\geq 1} \chi(\mathcal{M}^{s, B}_{K(2)}) x^{k (a+1)}y^{k a} = (1 + x^{a+1} y^{a})^{a+1}.
\end{align}
So the generating series are just polynomials in $x, y$, except for $B_{1,1}$.
\end{rmk}
\subsection{Sorting trees}\label{naiveTrees} Assumption 1 and 2 say that the sorting diagrams $\mathfrak{S}^i$ are well defined and stabilize. Assumption 3 says that it is easy to compose operators with the same slope in the stable sorting diagram $\mathfrak{S}^{\infty}$. We now spell out a further condition, which allows us to associate a tree $\overline{\Gamma}_{\sigma}$ with each element $\sigma \in \mathfrak{S}^i$. Again we will see that this assumption holds automatically when we work over the more general base rings of the next section. Our definition follows that of the tree underlying the scattering diagrams of \cite{gps} Section 1.4.  

Suppose that $\sigma \in \mathfrak{S}^i$ arises as the commutator of $\sigma_1, \sigma_2 \in \mathfrak{S}^{i-1}$. We define
\begin{equation}
\operatorname{Parents}(\sigma) = \{\sigma_1, \sigma_2\}.
\end{equation}
We then have the recursive functions
\begin{equation}
\operatorname{Ancestors}(\sigma) = \{\sigma\} \cup \bigcup_{\sigma' \in \operatorname{Parents}(\sigma)}\operatorname{Ancestors}(\sigma')
\end{equation}
and
\begin{equation}
\operatorname{Leaves}(\sigma) = \{\sigma' \in \operatorname{Ancestors}(\sigma) : \sigma' \in \mathfrak{S}^0\}.
\end{equation}

\textbf{Assumption 4.} If $\sigma' \in \operatorname{Ancestors}(\sigma)\setminus \{\sigma\}$, then $\sigma'$ is parent to a unique element of $\operatorname{Ancestors}(\sigma)$. We denote this by $\operatorname{Child}(\sigma')$.\\
\textbf{Example.} Going back to our examples in the previous section, we see that again this assumption holds in the first two cases of $Q_1, Q_2$, and fails for the subquiver $Q_3 \subset \widetilde{K}(3)$ that we considered. This is because $\mathfrak{S}^7$ contains the element $T_{\eta} = T_{i_1 +2 i_2 + i_3 + i_5 + i_6 + i_7}$. Then clearly $T_{i_2}$ must be parent to two different ancestors of $T_{\eta}$ (a little computation shows that these are in fact $T_{i_2 + i_5}, T_{i_2 + i_7}$).   

As before, Assumption 4 holds for subquivers of $\widehat{K}(2)$ without the need to pass to more general base rings.
\begin{lem}\label{special2} Let $Q \subset \widehat{K}(2)$ be a finite subquiver and let $\sigma$ be an element of some sorting diagram $\mathfrak{S}^i$. If $\sigma' \in \operatorname{Ancestors}(\sigma)\setminus \{\sigma\}$, then $\sigma'$ is parent to a unique element of $\operatorname{Ancestors}(\sigma)$. 
\end{lem}
\begin{proof} If there is $\sigma' \in \operatorname{Ancestors}(\sigma)$ violating Assumption 4, then we can actually choose $\sigma' = T_{i_p}$ for some $i_p \in Q_0$. It follows that for some $i > 0$, $\mathfrak{S}^i$ contains an element of the form $T_{d_A + d_B}$ for some dimension vectors of the form $d_A = \sum_{i' \in A} i'$, $d_B = \sum_{i'' \in B} i''$  for some distinct $A, B \subset Q_0$ with $i_p \in A \cap B$. As before this contradicts the dimension formula \eqref{dim}.
\end{proof}
Under Assumption 4, we proceed to construct the tree $\overline{\Gamma}_{\sigma}$. Both the vertices and edges of $\overline{\Gamma}_{\sigma}$ are parametrized by ancestors of $\sigma$:
\begin{align}
\overline{\Gamma}^{[0]}_{\sigma} &= \{V_{\sigma'} : \sigma' \in \operatorname{Ancestors}(\sigma)\},\\
\overline{\Gamma}^{[1]}_{\sigma} &= \{E_{\sigma'} : \sigma' \in \operatorname{Ancestors}(\sigma)\setminus \{\sigma\}\}.
\end{align}
Then for $\sigma' \in \operatorname{Ancestors}(\sigma)\setminus \{\sigma\}$ the vertices of $E_{\sigma'}$ are $\{V_{\sigma'}, V_{\operatorname{Child}(\sigma')}\}$.\\
\textbf{Example.} Going back to our two examples for $\widehat{K}(2)$ in the previous section, the tree for $\overline{\Gamma}_{T_{i_1+i_2+i_3}}$ in $Q_1$ is given by 
\begin{center}
\centerline{
\xymatrix{
{T_{i_2}}\ar[dr]&                      &{T_{i_3}}\ar[dl] &\\
                &{T_{i_2 + i_3}}\ar[dr]&                 &{T_{i_1}}\ar[dl]\\
                &                      &{T_{i_1+i_2+i_3}} &
}}
\end{center}
while the tree $\overline{\Gamma}_{T_{i_1+i_2+i_3 + i_4 + i_5}}$ in $Q_2$ is
\begin{center}
\centerline{
\xymatrix{
                                            &T_{i_2}\ar[d]                &T_{i_4}\ar[dl]\\    
                T_{i_1}\ar[d]               &{T_{i_2 + i_4}}\ar[dl]       &T_{i_3}\ar[d] & T_{i_5}\ar[dl]\\
                {T_{i_1 + i_2 + i_4}}\ar[dr]&                             &{T_{i_3 + i_5}}\ar[dl]\\
                                            &{T_{i_1+i_2+i_3 + i_4 +i_5}} &
}}
\end{center}

We also define a related \emph{unbounded} tree $\Gamma_{\sigma}$. We set
\begin{align}
\Gamma^{[0]}_{\sigma} &= \{V_{\sigma'} : \sigma' \in \operatorname{Ancestors}(\sigma)\setminus \mathfrak{S}^0\},\\
\Gamma^{[1]}_{\sigma} &= \{E_{\sigma'} : \sigma' \in \operatorname{Ancestors}(\sigma)\}.
\end{align}
For $\sigma' \in \operatorname{Ancestors}(\sigma)\setminus(\{\sigma\}\cup\mathfrak{S}^0)$, the vertices of $E_{\sigma'}$ are defined as $\{V_{\sigma'}, V_{\operatorname{Child}(\sigma')}\}$ as before. However, for $\sigma' \in \operatorname{Leaves}(\sigma)$ we define $E_{\sigma'}$ to be an unbounded edge with the single vertex $V_{\operatorname{Child}(\sigma')}$. Similarly we define $E_{\sigma}$ to be unbounded, with the single vertex $V_{\sigma}$. 

We define the \emph{weight} on edges $w_{\Gamma_{\sigma}}\!: \Gamma^{[1]}_{\sigma} \to \N_{> 0}$ as follows. We know that $\sigma' \in \operatorname{Ancestors}(\sigma)$ is a group element of the form $T_d$ for some $d \in \N Q_0$. Recall that we defined a reduction $\overline{d} \in \N K(m)$. We set
\begin{equation}\label{weight}
w_{\Gamma_{\sigma}}(E_{\sigma'}) = \ind(\overline{d}).
\end{equation}
\begin{rmk} The reason that we further reduce to $\overline{d}$ is that we will be interested in constucting \emph{plane tropical curves} from $\Gamma_{\sigma}$. Their integral structure is modelled on the rank 2 lattice $\Z K(m)_0$ rather then the higher rank lattice $\Z Q_0$.
\end{rmk} 
\subsection{More general base rings}\label{baseRings} In this section we apply the factorization-deformation technique developed in \cite{gps} Section 1.4. The main advantage is that our Assumptions 1, 2 and 3 will hold automatically in this context. On the other hand Assumption 4 does not hold in general, but we will see that composing operators in the stable sorting diagram is still reasonably simple.  Geometrically, in section \ref{curveConstruction} this will give rise to correction terms coming from \emph{disconnected} curves.

We introduce auxiliary variables $t_1, \dots, t_s, t_{s+1}, \dots, t_S$, and redefine the Kontse\-vich-Soibelman operators as elements of the group 
\begin{equation}
\Aut_{\C[[t_{\bullet}]]}\C[x_1,x^{-1}_1 \dots, x_S, x^{-1}_S][[t_{\bullet}]],
\end{equation} 
given by 
\begin{equation}
T_{d}(x_i) = x_i (1 + t^d x^d)^{\bra d, i\ket},
\end{equation}
for $d \in \N Q_0$. Fix an integer $k \geq 1$. We will work modulo the ideal
\begin{equation}
(t^{k+1}_1, \dots, t^{k+1}_{s}, t^{k+1}_{s+1}, \dots, t^{k+1}_{s + S}).
\end{equation} 
The full information of the Kontsevich-Soibelman operators is recovered in the limit $k \to \infty$. To make this precise define the ring
\begin{equation}
R_k = \C[[t_1, \dots, t_{s}, t_{s+1}, \dots, t_{s + S}]]/(t^{k+1}_1, \dots, t^{k+1}_{s}, t^{k+1}_{s+1}, \dots, t^{k+1}_{s + S}).
\end{equation}
We describe how to pass to a version of the fundamental product \eqref{product} which plays the same role as the standard scattering diagrams of \cite{gps} Definition 1.10. Consider the ring
\begin{equation}
\widetilde{R}_k = \C[\{u_{ij}, 1 \leq i \leq s + S, 1 \leq j \leq k\}]/(u^2_{ij}, 1 \leq i \leq s + S, 1 \leq j \leq k).
\end{equation}
There is an inclusion $R_k \hookrightarrow \widetilde{R}_k$ induced by 
\begin{equation}
t_i \mapsto \sum^k_{j = 1} u_{ij}.
\end{equation}
We can factor each of the operators $T_{i}$ in \eqref{product} over $\widetilde{R}_k$. First we have the identity in $R_{k}$,
\begin{equation}
\log(1 + t_i x_i) = \sum^k_{j = 1} \frac{(-1)^{j-1}}{j}t^j_i x^j_i.
\end{equation}
Now in $\widetilde{R}_k$,
\begin{equation}
t^j_i = \sum_{J \subset \{1, \dots, k\}, \# J = j} j! \prod_{l \in J} u_{il}.
\end{equation} 
Therefore 
\begin{equation*}
\log(1 + t_i x_i) = \sum^k_{j = 1} \sum_{J \subset \{1, \dots, k\}, \# J = j} (-1)^{j-1}(j-1)!\prod_{l \in J} u_{il}\,x^j_i, 
\end{equation*}
and since the variables $u_{il}$ are $2$-nilpotent, 
\begin{align*}
1+ t_i x_i &= 1 + \big(\sum^k_{l = 1} u_{il}\big) x_i\\
&= \prod^k_{j=1} \prod_{J \subset \{1, \dots, k\}, \# J = j} \big(1 + (-1)^{j-1}(j-1)! \prod_{l \in J} u_{il}\,x^j_i\big). 
\end{align*} 
This leads to the factorisation
\begin{equation}
T_i \equiv \prod_{J \subset \{1,\dots, k\}}T_{i,J} \mod (t^{k+1}_1, \dots, t^{k+1}_{s+S}),
\end{equation}
where the operators $T_{i, J}$ act by 
\begin{equation}\label{factoredOps}
T_{i, J}(x_j) = x_j (1 + (-1)^{\# J-1}(\# J-1)! \prod_{l \in J} u_{il}\,x^{\# J}_i)^{\bra i, j \ket}.
\end{equation}
Notice that $[T_{i, J}, T_{i, J'}] = 0$, so $\prod_J T_{i, J}$ is well defined. 

For any subset 
\[I \subset \{1, \dots, s + S\} \times \{1, \dots, k\}\]
we introduce the notation 
\begin{equation}
u_{I} = \prod_{(i,j) \in I} u_{ij}.
\end{equation} 
The following computation should be compared to \cite{gps} Lemma 1.9.
\begin{lem}\label{Commutator} Let $d_1, d_2$ be two primitive dimension vectors. Consider two operators $A_1, A_2$ acting by
\begin{align}
\nonumber A_1(x_j) &= x_j(1 + c_1 u_{I_1} x^{r_1 d_1})^{\bra d_1, j\ket},\\ 
A_2(x_j) &= x_j(1 + c_2 u_{I_2} x^{r_2 d_2})^{\bra d_2, j\ket}
\end{align} 
for some $c_i \in \C$, $I_i \subset \{1, \dots, s + S\} \times \{1, \dots, k\}$, $r_i \in \N_{>0}$, $i = 1, 2$. Then 
\begin{equation}
A^{-1}_2 \circ A_1 \circ A_2\circ A^{-1}_1 = B,
\end{equation}
where the operator $B$ acts by
\begin{equation}\label{operatorB}
B(x_j) = x_j(1 + c_1 c_2 \ind(r_1 d_1+r_2 d_2)\bra d_1, d_2\ket u_{I_1 \cup I_2}x^{r_1d_1 + r_2d_2})^{\bra \frac{r_1d_1 + r_2d_2}{\ind(r_1 d_1+r_2 d_2)}, j\ket}.
\end{equation}
In particular, if $I_1 \cap I_2 \neq \emptyset$, then the operators $A_1, A_2$ commute.
\end{lem}
\begin{proof} It is convenient to write $A_1, A_2$ as exponentials of derivations of the noncommutative Poisson algebra, 
\begin{equation}
A_1 = \exp\left(\{\frac{c_1}{r_1} u_{I_1} x^{r_1 d_1}, \cdot\,\}\right),\,\,\, A_2 = \exp\left(\{\frac{c_2}{r_2} u_{I_2} x^{r_2 d_2},\,\cdot\}\right).
\end{equation}
However since for $\xi, \eta$ in the Poisson algebra we have
\begin{equation}
[\{\xi,\, \cdot\}, \{\eta,\,\cdot\}] = \{\{\xi, \eta\},\,\cdot\}
\end{equation}
we will be sloppy and identify $\{\xi,\,\cdot\}$ with $\xi$ in the following. We compute (using the Baker-Campbell-Hausdorff formula and nilpotency)
\begin{align}
A_1 \circ A_2 &= \exp\left(\frac{c_1}{r_1} u_{I_1}x^{r_1 d_1} + \frac{c_2}{r_2} u_{I_2}x^{r_2 d_2} + \frac{1}{2}\bra d_1, d_2 \ket c_1 c_2 u_{I_1 \cup I_2} x^{r_1d_1 + r_2d_2}\right),\label{compose}\\
\nonumber A^{-1}_2 \circ A_1 \circ A_2 &= \exp\left(-\frac{c_2}{r_2} u_{I_2}x^{r_2 d_2} + \frac{c_1}{r_1} u_{I_1}x^{r_1 d_1} + \frac{c_2}{r_2} u_{I_2}x^{r_2 d_2}\right.\\
\nonumber &\left. + \frac{1}{2}\bra d_1, d_2 \ket c_1 c_2 u_{I_1 \cup I_2} x^{r_1d_1 + r_2d_2} - \frac{1}{2}\bra d_2, d_1 \ket c_1 c_2 u_{I_1 \cup I_2} x^{r_1d_1 + r_2d_2}\right)\\
&= \exp\left(\frac{c_1}{r_1} u_{I_1}x^{r_1 d_1} + \bra d_1, d_2 \ket c_1 c_2 u_{I_1 \cup I_2} x^{r_1d_1 + r_2d_2}\right), 
\end{align}
therefore
\begin{align}
\nonumber A^{-1}_2 \circ A_1 \circ A_2 \circ A^{-1}_1 &= \exp\left(\frac{c_1}{r_1} u_{I_1}x^{r_1 d_1} + \bra d_1, d_2 \ket c_1 c_2 u_{I_1 \cup I_2} x^{r_1d_1 + r_2d_2}\right.\\ 
\nonumber &\left.\,\,\,\,\,\,\,\,\,\,\,\,\,\,\,\,\,\,\,-\frac{c_1}{r_1} u_{I_1}x^{r_1 d_1} \right)\\
&= \exp\left(\bra d_1, d_2 \ket c_1 c_2 u_{I_1 \cup I_2} x^{r_1d_1 + r_2d_2}\right),
\end{align}
which in turn is identified with $\exp\left(\{\bra d_1, d_2 \ket c_1 c_2 u_{I_1 \cup I_2} x^{r_1d_1 + r_2d_2},\,\cdot\}\right)$, acting as in \eqref{operatorB}.
\end{proof}
\begin{rmk} The formula \eqref{compose} for the composition of operators which we obtained in the course of the proof will play a very important role in the following. 
\end{rmk}

We can now define the notions of \emph{sorting diagrams} $\mathfrak{S}^i_k$ over $\widetilde{R}_k$. These are ordered sequences of Poisson automorphisms $\sigma^i_j$ over $\widetilde{R}_k$, 
\begin{equation}
\mathfrak{S}^i_k = (\sigma^i_1, \dots, \sigma^i_{\ell_i}).
\end{equation}
(of course $\sigma^i_j, \ell_i$ also depend on $k$, but we omit this in the notation for brevity). We set 
\begin{equation}
\mathfrak{S}^0_k := ((T_{i_1, J_1})_{J_1 \subset \{1,\dots, k\}}, \ldots, (T_{i_S, J_S})_{J_S \subset \{1,\dots, k\}}\},
\end{equation}
where for each subsequence $(T_{i_S, J_S})_{J_S \subset \{1,\dots, k\}}$ we choose the \emph{lexicographic} order induced by subsets of $\{1, \dots, k\}$. Notice that for all elements of $\mathfrak{S}^0$ we have a well defined notion of slope: since $\sigma^0_j = T_{i_j, J}$ for some vertex $i_j$ and subset $J \subset \{1, \dots, k\}$, we set 
\begin{equation}
\mu(\sigma^0_j) = \mu((\# J)! i_j) = \mu(i_j).
\end{equation}
We then define the $\mathfrak{S}^i$ for $i > 0$ inductively precisely as in section \ref{naiveDiagrams}, with the only caveat that in commuting an element $\sigma^i_q$ past $\sigma^i_{q'}$ with $\mu(\sigma^i_{q}) < \mu(\sigma^i_{q'})$ we use Lemma \ref{Commutator} in place of Assumption 1. In particular by induction $\sigma^{i+1}_j \in \mathfrak{S}^{i+1}$ is an operator of the form 
\begin{equation}\label{sigma}
\sigma^{i+1}_j(x_p) = x_p(1+ c\,u_{I} x^{r d})^{\bra d, p\ket}
\end{equation}
for some $c \in \C, I \subset \{1, \dots, s + S\}\times \{1, \dots, k\}$ and primitive $d \in \N Q_0$, and we can define the slope 
\begin{equation}
\mu(\sigma^{i+1}_j) = \mu(r d) = \mu(d).
\end{equation} 
So we have well defined sorting diagrams $\mathfrak{S}^i_k$ for $i > 0$. Notice that since $u_{I_1 \cup I_2} = 0$ if $I_1 \cap I_2 \neq \emptyset$, the $\mathfrak{S}^{i}_{k}$ stabilise for $i > (s + S)k$, i.e. our Assumption 2 holds. 

As in section \ref{naiveTrees} we define the (bounded and unbounded) \emph{sorting trees} $\overline{\Gamma}_{\sigma}$ and $\Gamma_{\sigma}$ for $\sigma \in \mathfrak{S}^i$. First, the recursive functions $\operatorname{Ancestors}(\sigma)$, $\operatorname{Parents}(\sigma)$, $\operatorname{Leaves}(\sigma)$ are defined exactly as before. This extends immediately the definition of $\overline{\Gamma}_{\sigma}$ to the present case of diagrams over $\widetilde{R}_k$, provided we can show that our Assumption 4 holds. Namely, if $\sigma' \in \operatorname{Ancestors}(\sigma)$, we must show that it is parent to a unique element of $\operatorname{Ancestors}(\sigma)$. By \eqref{sigma} we know that $\sigma'$ is an operator of the form $x_p \mapsto (1 + c u_I x^{r d})^{\bra d, p\ket}$ for some nonempty set $I$. By Lemma \ref{Commutator} and induction, all its descendents must be operators of the form $x_p \mapsto (1 + c u_{I'} x^{r' d'})^{\bra d', p\ket}$ where $I \subset I'$. Applying again Lemma \ref{Commutator} we see that two operators of this form commute. Thefore at most one descendent of $\sigma'$ appears in $\operatorname{Ancestors}(\sigma)$. The unbounded tree $\Gamma_{\sigma}$ is then obtained from $\overline{\Gamma}_{\sigma}$ exactly as before. 

Finally we define the \emph{weight} over $\widetilde{R}_k$, $w_{\Gamma_{\sigma}}\!: \Gamma_{\sigma} \to \N_{> 0}$. By \eqref{sigma} we can write $\sigma' \in \operatorname{Ancestors}(\sigma)$ uniquely in the form $x_p \mapsto (1 + c u_I x^{r d})^{\bra d, p\ket}$. Then we set
\begin{equation}\label{weightGeneral}
w_{\Gamma_{\sigma}}(E_{\sigma'}) = \ind(\overline{r d}) = r \ind(\overline{d}). 
\end{equation}

\subsection{Tropical curves and counts}\label{preTropical} In this section we recall the notions of rational tropical curves in $\R^2$ and of their counting invariant. We follow \cite{gps} Sections 2.1 and 2.3. 

Let $\Gamma$ be a weighted, unbounded tree with only trivalent vertices. We have a weight 
\begin{equation}
w_{\Gamma}\!: \Gamma^{[1]} \to \N_{> 0}
\end{equation}
and a distinguished subset of noncompact edges $\Gamma^{[1]}_{\infty} \subset \Gamma^{[1]}$ (which as usual we call unbounded edges). A \emph{parametrized rational tropical curve} in $\R^2$ is a proper map $h\!: \Gamma \to \R^2$ such that:
\begin{enumerate}
\item[$\bullet$] for every $E \in \Gamma^{[1]}$, the restriction $h|_{E}\!: E \to \R^2$ is an embedding with image $h(E)$ contained in an affine line of rational slope;
\item[$\bullet$] for every $V \in \Gamma^{[0]}$, if $E_i, i = 1, 2, 3$ are the edges adjecent to $V$ and $m_i, i = 1, 2, 3$ is the primitive integral vector emanating from $h(V)$ in the direction of $h(E_i)$, we have the balancing condition
\begin{equation}\label{balance}
w_{\Gamma}(E_1) m_1 + w_{\Gamma}(E_2) m_2 + w_{\Gamma}(E_3) m_3 = 0.
\end{equation}
\end{enumerate}
Two parametrized rational tropical curves $h\!: \Gamma \to \R^2, h'\!: \Gamma' \to \R^2,$ are \emph{equivalent} if there is a homeomorphism $\Phi\!: \Gamma \to \Gamma'$, respecting the weights of the edges, such that $h'\circ\Phi = h$. A \emph{rational tropical curve} is an equivalence class of parametrized rational tropical curves.

Following the notation of balancing condition, we define the \emph{multiplicity of a vertex} $V$ as
\begin{equation}\label{multipVert}
\mult_V(h) = w_{\Gamma}(E_i)w_{\Gamma}(E_j)|m_i\wed m_j|
\end{equation}  
for $i \neq j$. This gives a good definition by the balancing condition. The \emph{multiplicity of a tropical curve} $h$ is then defined as
\begin{equation}\label{multipCurve}
\mult(h) = \prod_V \mult_V(h).
\end{equation}

Write $m_1, \dots, m_n$ for primitive vectors of $\R^2$ (not necessarily distinct), and $\m$ for their $n$-tuple. Choose generic vectors $m_{ij} \in \R^2$ for $1 \leq i \leq n$,  $1 \leq j \leq l_i$, and form the lines
\begin{equation}
\mathfrak{d}_{ij} = m_{ij} + \R m_i \subset \R^2.
\end{equation} 
Let $\w_i = (w_{i1}, \dots, w_{il_i})$, $1 \leq i \leq n$ be \emph{weight vectors} with 
\begin{equation}
0 < w_{i1} \leq w_{i2} \leq \dots \leq w_{i l_i}.
\end{equation}
The weight vector $\w_i$ has \emph{length} $l_i$ and \emph{size} $|\w_i| = \sum w_{ij}$. We also need the notion of the \emph{automorphism group} of a weight vector: this is the subgroup $\Aut(\w_i)$ of the permutation group $\Sigma_{l_i}$ stabilizing the vector $(w_{i1}, \dots, w_{il_i})$. We will write $\w = (\w_1, \dots, \w_n)$ and set $\Aut(\w) = \prod^n_{i = 1} \Aut(\w_i)$. We will also use the notation $m_{\out} = \sum_i |\w_i|m_i$.  

Consider the (finite) set $\mathcal{S}(\w)$ of tropical curves $h\!: \Gamma \to \R^2$ which satisfy the following constraints:
\begin{enumerate}
\item[$\bullet$] the unbounded edges of $\Gamma$ are 
\begin{equation}
\Gamma^{[1]}_{\infty} = \{E_{ij}, 1 \leq i \leq n, 1\leq j \leq l_i\}\cup\{E_{\out}\},
\end{equation}
and $h(E_{ij}) \subset \mathfrak{d}_{ij}$, with $-m_i$ pointing in the unbounded direction of $h(E_{ij})$, and $m_{\out}$ pointing in the unbounded direction of $h(E_{\out})$;
\item[$\bullet$] $w_{\Gamma}(E_{ij}) = w_{ij}$.  
\end{enumerate}
\begin{thm}[\cite{gps} Proposition 2.7] The number of elements of $\mathcal{S}(\emph{\w})$, counted with the multiplicity of \eqref{multipCurve}, is independent of the generic choice of lines (i.e. of the vectors $m_{ij}$). We denote this number by $N^{\trop}_{\emph{\m}}(\emph{\w})$.
\end{thm}

The Gromov-Witten invariants which appear in the GW/Kronecker correspondence \eqref{KronGW} arise from tropical counts. Fix an ordered partition $P = (P_1, P_2)$ (we wrote $(P_a, P_b)$ for this in the Introduction). Choose $m_1 = (1, 0), m_2 = (0,1)$. For this standard choice we omit $\m$ from the notation. A weight vector $\w = (\w_1, \w_2)$ has the same \emph{type} as $P$ if $|P_i| = |\w_i|$ for $i = 1, 2$. In this case we write $\w \sim P$.  Let us write $p_{ij}$ for the (ordered) parts of $P_i$, and $I_{\bullet}$ for a partition of the sets $\{1, \dots, l_i\}$:
\begin{equation}
I_1 \cup \dots \cup I_{l_{i}} = \{1, 2, \dots, l_i\}.
\end{equation} 
We call $I_{\bullet}$ a set partition of $\w_i$, and say it is \emph{compatible} with $P_i$ if 
\begin{equation}
p_{ij} = \sum_{r \in I_j} w_{ir}.
\end{equation}
For an integer $r > 0$, we set $R_r = \frac{(-1)^{r-1}}{r^2}$, and we define some coefficients
\begin{equation}
R_{P_i | \w_i} = \sum_{I_{\bullet}} \prod^{l_i}_{j = 1} R_{w_{ij}},
\end{equation}
where the sum is over all set partitions of $\w_i$ which are compatible with $P_i$. Set $R_{P | \w} = \prod_{i = 1, 2} R_{P_i |\w_i}$.
\begin{thm}[\cite{gps} Theorem 3.4 and Proposition 5.3] 
\begin{equation}\label{tropicalThm}
N(P) = \sum_{\emph{\w}\sim P}\frac{R_{P | \emph{\w}}}{|\Aut(\emph{\w})|}N^{\trop}(\emph{\w}).
\end{equation}
\end{thm}
As we mentioned, one of the heuristic motivations for the present work is to compare the formula \eqref{tropicalThm} to Weist's result \eqref{weist}. In other words we would like to regard a graded partition $P$ as analogue to a dimension vector $\overline{d} \in \N K(m)_0$, and a weight vector $\w$ which has the same type as $P$ as an analogue of a dimension vector $d \in \N\widetilde{K}(m)_0$ which is compatible with $d$. We achieve this at least in part in section \ref{main}.
\subsection{Tropical curves and counts from $\widetilde{K}(m)$}\label{curveConstruction} Fix $Q \subset \widetilde{K}(m)$. Let $\sigma \in \mathfrak{S}^{\infty}$ be an element of the stable sorting diagram. To this we will associate a rational tropical curve $h_{\sigma}\!: \Gamma_{\sigma} \to \R^2$, together with a dimension vector $d_{\out}(h_{\sigma}) \in \N Q_0$. 

\begin{figure}[ht]
\centerline{\includegraphics[scale=.7]{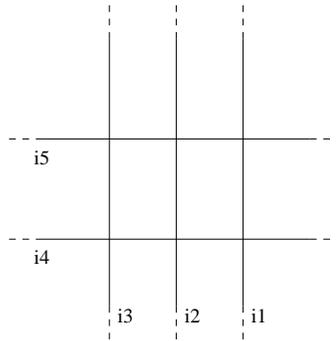}}
\caption{Line arrangement.}
\label{lines}
\end{figure}

We start in the simplified situation of sections \ref{naiveDiagrams} and \ref{naiveTrees}, where the naive definition of sorting diagrams and trees apply. Fix $s$ vertical lines, $S$ horizontal lines in $\R^2$ generically. We label the vertical lines $\mathfrak{d}_{i_j}$ with $i_1, \dots, i_s$ starting from the \emph{rightmost} line. Similarly we label the horizontal lines $\mathfrak{d}_{i_j}$ with $i_{s+1}, \dots, i_{s+S}$ starting from the \emph{lowest} line. In other words, the lines $\mathfrak{d}_{i_j}$ are labelled with $i_1, \dots, i_{s+S}$ in \emph{clockwise} order starting from the rightmost vertical line.\\
\textbf{Example} The line arrangement for the quiver $Q_2$ of section \ref{naiveDiagrams} is shown in Figure \ref{lines}.

\begin{figure}[ht]
\centerline{\includegraphics[scale=.7]{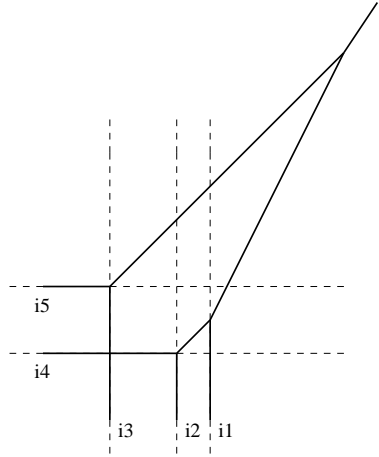}}
\caption{The tropical curve for $T_{i_1 + i_2 +i_3 + i_4 + i_5}$.}
\label{lines23}
\end{figure}

\begin{figure}[ht]
\centerline{\includegraphics[scale=.7]{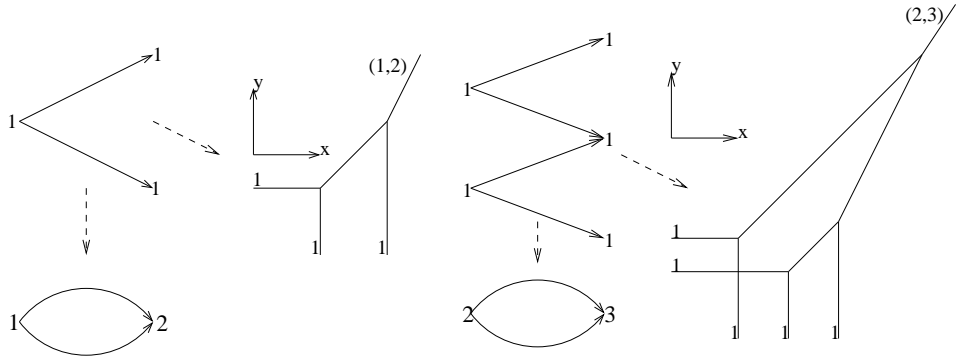}}
\caption{From quivers to tropical curves for $K(2)$.}
\label{transform}
\end{figure}

Pick $\sigma \in \mathfrak{S}^{\infty}$. It will appear for the first time in the sequence of diagrams $\mathfrak{S}^i$ for some finite $i \geq 0$. If $i = 0$ we are in a degenerate case, $\sigma = T_{i_p}$ for some $p$ and we just choose $h_{\sigma}$ to be the corresponding (vertical or horizontal) line. In this case we also set $d_{\out}(h) = i_p$. If $\sigma$ first appears in $\mathfrak{S}^1$ and $\operatorname{Parents}(\sigma) = \{T_{i_p},  T_{i_q}\}$ with $p < q$ then $1 \leq p \leq s$ and $s + 1 \leq q \leq S$. We define $h_{\sigma}$ as the unique tropical curve $h_{\sigma}\!:\Gamma_{\sigma} \to \R^2$ with unbounded vertical edge $\mathfrak{d}_{i_{p}}$ and unbounded horizontal edge $\mathfrak{d}_{i_{q}}$. We also set $d_{\out}(h) = i_p + i_q$. Suppose now $\sigma$ first appears in $\mathfrak{S}^i$ with $i > 1$. We must have $\operatorname{Parents}(\sigma) = \{\sigma_1, \sigma_2\}$ with $\sigma_1, \sigma_2 \in \mathfrak{S}^{i-1}$ and $\mu(\sigma_1) < \mu(\sigma_2)$. By induction we have well defined rational tropical curves $h_{\sigma_1}\!: \Gamma_{\sigma_1} \to \R^2$ and $h_{\sigma_2}\!: \Gamma_{\sigma_2} \to \R^2$. Notice that by construction the slope inequality $\mu(\sigma_1) < \mu(\sigma_2)$ (using the slope for quiver dimension vectors) implies the \emph{opposite} inequality for the slopes of the outgoing unbounded edges of the tropical curves. Namely, $\operatorname{slope}(h(E_{\sigma_1})) > \operatorname{slope}(h(E_{\sigma_2}))$ as rays in $\R^2$. But notice also that by our choice of labels for the quiver $Q$, the set $\operatorname{Leaves}(\sigma_1)$ preceeds the set $\operatorname{Leaves}(\sigma_2)$ in the lexicographic order and so by our choice of labels for the lines $\mathfrak{d}_{ij}$ the set $h(\operatorname{Leaves}(\sigma_1))$ preceeds $h(\operatorname{Leaves}(\sigma_2))$ in the clockwise order in $\R^2$. This implies that the ray $h(E_{\sigma_1})$ emanates from a point which lies below the ray $h(E_{\sigma_2})$. Therefore the two rays must intersect in $\R^2$. We then use the balancing condition given by the weights \eqref{weight} to construct $h_{\sigma}$ inductively as a map from $\Gamma_{\sigma}$. We can also define $d_{\out}(h_{\sigma})$ inductively as $d_{\out}(h_{\sigma_1}) + d_{\out}(h_{\sigma_2})$. For a curve $h$ corresponding to some $\sigma \in \mathfrak{S}^{\infty}$, we will write $\operatorname{Legs}(h)$ for the set of lines $\mathfrak{d}_{i}$, $i \in Q_0$ appearing in $h$.\\
\textbf{Example.} Consider once again the examples of sections \ref{naiveDiagrams}, \ref{naiveTrees}. For $Q_2 \subset \widehat{K}(2)$ the tree $\Gamma_{T_{i_1 + i_2 +i_3 + i_4 + i_5}}$ maps to the curve in Figure \ref{lines23}. On the other hand we can identify $Q_1$ with the subquiver of $Q_2$ spanned by $i_1, i_2, i_4$, and $\Gamma_{T_{i_1 + i_2 +i_4}}$ maps to the bottom subcurve with legs $i_1, i_2, i_4$. In general for $Q \subset \widehat{K}(2)$ we know by the proof of Lemma \ref{special1} that we can identify operators in $\mathfrak{S}^{\infty}$ with subquivers of $Q$, which then map to tropical curves by the construction above, see Figure \ref{transform} for a schematic picture.  

As for all tropical curves, we have the notion of multiplicity at a vertex $\mult_V h_{\sigma}$. We modify the notion of multiplicity using the quiver $Q$ as follows. A vertex $V \in h_{\sigma}$ corresponds to a pair of incoming dimension vectors $d_{V,1}, d_{V,2}$ with $\mu(d_{V,1}) < \mu(d_{V,2})$. We set 
\begin{equation}\label{QmultipVert}
\mult_{Q, V} h_{\sigma} = \bra d_{V,1}, d_{V,2}\ket.
\end{equation} 
The global multiplicity is
\begin{equation}\label{QmultipCurve}
\mult_{Q} h_{\sigma} = \prod_V \mult_{Q, V} h_{\sigma}.
\end{equation}
To compare with the usual notion of tropical multiplicity, notice that in fact
\begin{equation}\label{multipVert2}
\mult_V h_{\sigma} = \frac{1}{m}|\bra \overline{d}_{V,1}, \overline{d}_{V,2} \ket|
\end{equation}
(we still write $\bra \cdot,\,\cdot\ket$ for the product of the \emph{reduced} dimension vectors, computed on $K(m)$).\\ 
\textbf{Example.} In the example of subquivers $Q_1, Q_2 \subset \widehat{K}(2)$ of section \ref{naiveDiagrams} we have 
\[\mult_{Q_1}h_{T_{i_1 + i_2 + i_3}} = \mult h_{T_{i_1 + i_2 + i_3}} = 1,\]
and also 
\[\mult_{Q_2}h_{T_{i_1 + i_2 + i_3 + i_4 + i_5}} = \mult h_{T_{i_1 + i_2 + i_3+i_4+i_5}} = 1.\]
But we can compute
\[\mult_{Q_2} h_{T_{i_1 + i_2 + i_4 +i_5}} = 1, \mult h_{T_{i_1 + i_2 + i_4 +i_5}} = 2.\]

Finally, we will denote by $\mathcal{S}_Q$ the (finite) set of all rational tropical curves $h_{\sigma}$ that we constructed for $\sigma \in \mathfrak{S}^{\infty}$, and by $\mathcal{S}_Q(\mu)\subset \mathcal{S}_Q$ the subset of curves whose outgoing dimension vector has prescribed slope, namely $\mu(d_{\out}(h)) = \mu$.\\

\begin{figure}[ht]
\centerline{\includegraphics[scale=.7]{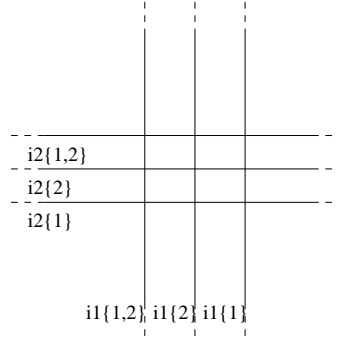}}
\caption{Line arrangement with $k = 2$.}
\label{lines2}
\end{figure}

We now move on the the general case, working over the base rings of section \ref{baseRings}. Recall in this case we have the additional parameter $k \geq 1$. We fix $(2^k-1)s$ vertical lines, $(2^k-1)S$ horizontal lines in $\R^2$ generically. We label the vertical lines $\mathfrak{d}_{i_j, I}$ with elements of $\{i_1, \dots, i_s\} \times \{I \subset \{1, \dots, k\}, I \neq \emptyset\}$ in lexicographic order, starting from the \emph{rightmost} line. Similarly we label the vertical lines $\mathfrak{d}_{i_j, J}$ with elements of $\{i_{s+1}, \dots, i_S\} \times \{J \subset \{1, \dots, k\}, J \neq \emptyset\}$ in lexicographic order, starting from the \emph{lowest} line. In other words, the set of all lines $\mathfrak{d}_{i_j, I}$ is labelled with $\{i_1, \dots, i_{s+S}\}\times\{I \subset \{1, \dots, k\}, I \neq \emptyset\}$ in \emph{clockwise} order starting from the rightmost vertical line, $\mathfrak{d}_{i_1, \{1\}}$ to the top horizontal line, $\mathfrak{d}_{i_{s+S}, \{1, \dots, k\}}$.\\
\textbf{Example.} The simplest case of a subquiver $i_2 \to i_1$ with $k = 2$ is shown in Figure \ref{lines2}.

Pick $\sigma \in \mathfrak{S}^{\infty}$. We want to construct a tropical curve $h_{\sigma}$ from $\sigma$. If $\sigma$ first appears in $\mathfrak{S}^0$ then $\sigma = T_{i_p, I}$ for some vertex $i_p$ and $I \subset\{1, \dots, k\}$. Then we are in a degenerate case and we just choose $h_{\sigma}$ to be the corresponding line $\mathfrak{d}_{i_p, I}$. If $\sigma$ first appears in $\mathfrak{S}^1$ and $\operatorname{Parents}(\sigma) = \{T_{i_p, I},  T_{i_q, J}\}$ with $p < q$ then $1 \leq p \leq s$ and $s + 1 \leq q \leq S$. We define $h_{\sigma}$ as the unique tropical curve $h_{\sigma}\!:\Gamma_{\sigma} \to \R^2$ with unbounded vertical edge $\mathfrak{d}_{i_{p}, I}$ and unbounded horizontal edge $\mathfrak{d}_{i_{q}, J}$. Then for $\sigma \in \mathfrak{S}^i$ with $i > 1$, we construct $h_{\sigma}$ inductively precisely as in the discussion above, using now the weight \eqref{weightGeneral} in the balancing condition \eqref{balance}. Similarly, we define the multiplicity of a vertex $V \in h_{\sigma}$ by \eqref{QmultipVert}, namely if the vertex arises from commuting $x_p \mapsto x_p(1 + c_1 x^{r_1 d_1})^{\bra d_1, p\ket}$ and $x_p \mapsto (1 + c_2 x^{r_2 d_2})^{\bra d_2, p\ket}$ its multiplicity is $\bra r_1 d_1, r_2 d_2\ket$. The global multiplicity of $h_{\sigma}$ is then given by \eqref{QmultipCurve}.

As we mentioned in section \ref{baseRings}, Assumption 4 in the definition of naive sorting diagrams does not hold in general for diagrams over $\widetilde{R}_k$. Correction terms will arise from a class of \emph{disconnected} tropical curves, i.e. maps from disconnected trees, which we now define. Fix any ordered $l$-uple of elements of $\mathfrak{S}^{\infty}_k$ \emph{with the same slope},
\begin{equation}
(\sigma_1, \dots, \sigma_l) \in \mathfrak{S}^{\infty}_k, \mu(\sigma_i) = \mu(\sigma_{i+1}),
\end{equation} 
such that $\sigma_i \prec \sigma_{i+1}$ in $\mathfrak{S}^{\infty}_k$, and the sets $I_1, \dots I_l$ underlying $\sigma_1, \dots, \sigma_l$ are pairwise disjoint. We define a tropical curve $h_{\sigma_1\cdots\sigma_l}$ simply as the union of the tropical curves $h_{\sigma_1}, \dots, h_{\sigma_l}$ (a map from the disjoint union $\cup^l_{i = 1} \Gamma_{\sigma_i}$). We still use the notation $\operatorname{Legs}(h_{\sigma_1\cdots\sigma_l})$ for the set of lines $\mathfrak{d}_{i, I}$ appearing in the image of $h$.

To each curve $h_{\sigma_1 \cdots \sigma_l}$ we attach inductively a \emph{weight function} $f_{\sigma_1 \cdots \sigma_l}$ as follows. We know from \eqref{sigma} that $\sigma$ acts by $x_p \mapsto x_p(1+ c u_I x^{r d})^{\bra d, i_p\ket}$ for some primitive $d$. We set 
\begin{equation}\label{elemWeight}
f_{\sigma} = \frac{c}{r} u_I x^{r d}.
\end{equation}
Suppose then inductively that $f_{\sigma_1\cdots \sigma_{l-1}} = \alpha u_{I} x^{r d}$ for some primitive $d$ and $\alpha \in \C$, and similarly $f_{\sigma_{l}} = \alpha' u_{I'} x^{r' d'}$. Then we set 
\begin{equation}\label{weightFunctDef}
f_{\sigma_1 \cdots \sigma_l} = \frac{1}{2}\bra r d, r' d'\ket \alpha\,\alpha' u_{I \cup I'} x^{r d + r'd'}.
\end{equation}
Representing $f_{\sigma_1 \cdots \sigma_l}$ in the form $c u_I x^d$, we set \begin{equation}
d_{\out}(h_{\sigma_1 \cdots \sigma_l}) = d.
\end{equation}
We extend the notion of multiplicity for $h_{\sigma_1 \cdots \sigma_l}$ inductively as follows:
\begin{equation}\label{QmultipCurve_disc}
\mult_Q(h_{\sigma_1 \cdots \sigma_l}) = \frac{1}{2}\mult_Q(h_{\sigma_1 \cdots \sigma_{l-1}})\bra d_{\out}(h_{\sigma_{1}\cdots\sigma_{l-1}}), d_{\out}(\sigma_l)\ket.
\end{equation}

We write $\mathcal{S}_{Q,k}$ for the set of all tropical curves $h_{\sigma_1 \cdots \sigma_l}$ (we do not fix $l$). For a fixed slope $\mu$, we write $\mathcal{S}_{Q, k}(\mu)$ for the subset of $\mathcal{S}_{Q, k}$ given by all curves $h_{\sigma_1 \cdots \sigma_l}$ with $\mu(d_{\out}(h_{\sigma_1 \cdots \sigma_l})) = \mu$.

Fix a curve $h_{\sigma_1 \cdots \sigma_l} \in \mathcal{S}_{Q, k}$. From this we find a unique weight vector 
\[\textbf{w} = (\w_1,\dots, \w_{s+S}),\]
\[\textbf{w}_q = (w_{q1}, \dots , w_{ql_q}),\]
with $1 \leq w_{p1} \leq \dots \leq w_{q l_q}$, and pairwise disjoint sets 
\[J_{qj} \subset \{1, \dots, k\}, q = 1, \dots, s+S;\, j = 1, \dots, l_q,\] with $\# J_{qj} = w_{qj},$ such that
\[\operatorname{Legs}(h_{\sigma_1 \cdots \sigma_l}) = \{\mathfrak{d}_{i_q, J_{qj}} |\, q = 1, \dots, s + S;\,j = 1, \dots, l_q\}.\]
Let us denote by $N^{\trop}_{Q, k}(\{J_{qj}\})$ the number of curves $h_{\sigma_1 \cdots \sigma_l}$ giving rise to the same sets $J_{qj}$, counted with the multiplicity \eqref{QmultipCurve_disc}. By the construction of $h_{\sigma_1 \cdots \sigma_l}$ in terms of sorting diagrams, it is clear that $N^{\trop}_{Q, k}(\{J_{qj}\})$ only depends on the vector $\w$, not the actual subsets $J_{qj}$. 

We denote this number by $N^{\trop}_{Q,k}(\w)$.

Notice that in fact $d_{\out} = \sum_i |\w_i|\,i$, so by abuse of notation we write $\mu(\w)$ and $\bra\cdot, \w\ket$ for the slope $\mu(d_{\out})$, respectively the linear form $\bra\cdot, d_{\out}\ket$. Similarly we will often write 
\begin{equation}
x^{\w} = x^{\sum_i |\w_i|\,i} = x^{|\w_1|}_{i_1} \cdots x^{|\w_{s+S}|}_{i_{s+S}}
\end{equation}
and 
\begin{equation}
[\,\w\,] = [\,\sum_i |\w_i|\,i\,],
\end{equation}
the equivalence class of the underlying dimension vector. For a fixed weight vector $\w$ we set 
\begin{equation}
R_{\textbf{w}} = \prod^{s+S}_{i = 1}\prod^{l_1}_{j = 1} \frac{(-1)^{w_{ij}-1}}{w^2_{ij}}
\end{equation}
and  
\begin{equation}
\overline{\w} = \overline{\sum_i |\w_i|\,i} = (\sum_{i > s}|\textbf{w}_i|, \sum_{i \leq s}|\textbf{w}_i|).
\end{equation}
\section{Main results}\label{main}
Throughout this section we will concentrate on $B$-framings for $K(m)$, and so on framings at sources on $\widetilde{K}(m)$. The situation for $F$-framings is completely analogous. We will use the notation $(f(x))[x^d]$ to denote the coefficient of $x^d$ in $f(x)$. We start in the simple situation described in sections \ref{naiveDiagrams}, \ref{naiveTrees}. According to Lemma \ref{special1} and Lemma \ref{special2}, our main example is a finite subquiver $Q \subset \widetilde{K}(2)$ (although soon we will need to restrict to representations of slope $\mu \neq \frac{1}{2}$ in order to have Assumption 3 in place). The following lemma simply summarizes our simplified construction in sections \ref{naiveDiagrams}, \ref{naiveTrees} and the first part of section \ref{curveConstruction}. 
\begin{lem} Suppose that Assumptions 1, 2, 4 hold. Then there is a bijective correspondence between operators $T_d$ with $\mu(d) = \mu$ appearing in the stable sorting diagram $\mathfrak{S}^{\infty}$ and tropical curves $h \in \mathcal{S}_Q(\mu)$. Moreover 
\begin{equation}
d = d_{\out}(h) = \sum_{\mathfrak{d}_i \in \operatorname{Legs}(h)} i.
\end{equation}
\end{lem}
\begin{lem}\label{simpleCompose} Suppose that Assumptions 1-4 hold. Let $T_{d_1} \prec \dots \prec T_{d_r}$ be the maximal sequence of operators with $\mu(d_i) = \mu$ appearing in $\mathfrak{S}^{\infty}$. Then the composition $T_{d_1}\circ\cdots\circ T_{d_r}$ acts by
\begin{equation}
x_p \mapsto x_p \prod_{h \in \mathcal{S}_Q(\mu)}(1 + x^{d_{\out}(h)})^{\bra d_{\out}(h), p\ket}.
\end{equation}
\end{lem}
\begin{proof} We know that $T_{d_i}$ corresponds to a curve $h \in \mathcal{S}_Q(\mu(d_i))$ and $d_i = d_{\out}(h)$, $T_{d_i}(x_p) = x_p (1 + x^{d_{\out}(h)})^{\bra d_{\out}(h), p\ket}$. This correspondence is bijective, and operators of the same slope $\mu(d_i)$ commute, and so compose simply as in the statement.
\end{proof}
\begin{cor} Let $Q \subset \widehat{K}(2)$ be a subquiver with $s + 1$ sinks, $s$ sources. Then for $1 \leq p \leq S$, $\mu \neq \frac{1}{2}$,
\begin{equation}\label{specialBackward}
\theta_{Q, \mu, i_{s + p}}(x) = \prod_{h \in \mathcal{S}_Q(\mu)} (1 + x^{d_{\out}(h)})^{\bra i_{p} - i_{p-1} + \dots \pm i_1, d_{\out}(h)\ket}.
\end{equation}
In other words
\begin{equation}
\chi(\mathcal{M}^{s, i_{s + p}}_Q(d)) = \prod_{h \in \mathcal{S}_Q(\mu)} (1 + x^{d_{\out}(h)})^{\bra i_{p} - i_{p-1} + \dots \pm i_1, d_{\out}(h)\ket}[x^d].
\end{equation}
\end{cor}
\begin{proof} First by Reineke's Theorem \ref{reinThm} (and our choice of a minimal labelling) we know
\begin{equation}
\theta_{Q, \mu}(x_{i_1}) = x_{i_1} (\theta_{Q, \mu, i_{s+1}}(x))^{-1}.
\end{equation}
By Lemma \ref{simpleCompose} we have
\begin{equation}
\theta_{Q, \mu}(x_{i_1}) = x_{i_1}\prod_{h \in \mathcal{S}_Q(\mu)}(1 + x^{d_{\out}(h)})^{\bra d_{\out}(h), i_1\ket},
\end{equation}
from which of course we find
\begin{equation}
\theta_{Q, \mu, i_{s+1}}(x) = \prod_{h \in \mathcal{S}_Q(\mu)}(1 + x^{d_{\out}(h)})^{\bra i_1, d_{\out}(h)\ket}.
\end{equation}
This establishes \eqref{specialBackward} for $p = 1$. For $p > 1$ we have again by Reineke's Theorem 
\begin{equation}
\theta_{Q, \mu}(x_{i_{p+1}}) = x_{i_{p+1}} (\theta_{Q, \mu, i_{s+p}}(x))^{-1}(\theta_{Q, \mu, i_{s+p+1}}(x))^{-1}.
\end{equation} 
Therefore
\begin{equation}
\theta_{Q, \mu, i_{s+p+1}}(x) = x_{i_{p+1}}(\theta_{Q, \mu}(x_{i_{p+1}}))^{-1}(\theta_{Q, \mu, i_{s+p}}(x))^{-1}. 
\end{equation}
By Lemma \ref{simpleCompose} we have
\begin{equation}
(\theta_{Q, \mu}(x_{i_{p+1}}))^{-1} = \frac{1}{x_{i_{p+1}}} \prod_{h \in \mathcal{S}_Q(\mu)}(1 + x^{d_{\out}(h)})^{\bra i_{p+1}, d_{\out}(h)\ket}.
\end{equation}
and assuming by induction that \eqref{specialBackward} holds for $(\theta_{Q, \mu, i_{s+p}}(x))^{-1}$ we find
\begin{align}
\nonumber\theta_{Q, \mu, i_{s+p+1}}(x) &= \prod_{h \in \mathcal{S}_Q(\mu)}(1 + x^{d_{\out}(h)})^{\bra i_{p+1}, d_{\out}(h)\ket} \prod_{h \in \mathcal{S}_Q(\mu)} (1 + x^{d_{\out}(h)})^{\bra - i_{p} + \dots \pm i_1, d_{\out}(h)\ket}\\
&= \prod_{h \in \mathcal{S}_Q(\mu)}(1 + x^{d_{\out}(h)})^{\bra i_{p+1} - i_p + \dots \pm i_1, d_{\out}(h)\ket}.
\end{align}
\end{proof}
\begin{cor}\label{contribSimple} Let $\overline{d}$ be a fixed dimension vector for $K(2)$. Choose $Q \subset \widehat{K}(2)$ with $s + 1$ sinks, $s$ sources for $s$ large enough (depending on $\overline{d}$). Then
\begin{equation}
\chi(\mathcal{M}^{s, B}_{K(2)}(\overline{d})) = \sum_{[d]\sim \overline d} \sum^s_{p = 1} \prod_{h \in \mathcal{S}_Q(\mu)} (1 + x^{d_{\out}(h)})^{\bra i_{p} - i_{p-1} + \dots \pm i_1, d_{\out}(h)\ket}[x^d],
\end{equation}
where the first sum is over all equivalence classes $[d]$ of dimension vectors supported on $Q$ and compatible with $\overline{d}$.
\end{cor}
\begin{proof} This is an immediate consequence of Weist's theorem in the form \eqref{weist} and the result we just proved, since the cardinality of the support of $d \in \N \widehat{K}(2)_0$ compatible with $\overline{d}$ is uniformly bounded (and so $d$ can be moved to an equivalent dimension vector in a large enough quiver $Q$ of the type we want). 
\end{proof}
\begin{figure}[ht]
\centerline{\includegraphics[scale=.7]{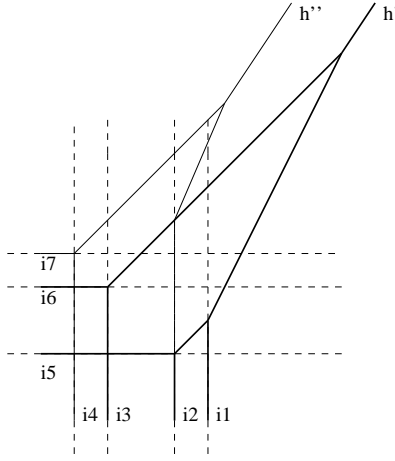}}
\caption{Curves contributing to $\chi(\mathcal{M}^{s, B}_{K(2)}(4, 6))$.}
\label{lines46}
\end{figure}
\begin{exa} By \eqref{genSeries} we know $B_{2, 3}(x, y) = (1 + x^2 y^3)^2$, so $\chi(\mathcal{M}^{s, B}_{K(2)})(4, 6) = 1$. This is witnessed already by a subquiver $Q \subset \widehat{K}(2)$ with $s = 4$, $S = 3$. A lengthy but elementary computation shows that $\mathcal{S}_Q(\frac{2}{5})$ contains only two curves $h', h''$ with $\operatorname{Legs}(h') = \{\mathfrak{d}_{i_1}, \mathfrak{d}_{i_2}, \mathfrak{d}_{i_3}, \mathfrak{d}_{i_5}, \mathfrak{d}_{i_6}\}$ and $\operatorname{Legs}(h'') = \{\mathfrak{d}_{i_2}, \mathfrak{d}_{i_3}, \mathfrak{d}_{i_4}, \mathfrak{d}_{i_6}, \mathfrak{d}_{i_7}\}$ (see Figure \ref{lines46}). So $d_{\out}(h') = i_1 + i_2 + i_3 + i_5 + i_6$ and $d_{\out}(h'') = i_2 + i_3 + i_4 + i_6 + i_7$, and 
\begin{equation}
\prod_{h \in \mathcal{S}_Q(\frac{2}{5})} (1 + x^{d_{\out}(h)})^{\bra i_2 - i_1, d_{\out}(h)\ket} = (1 + x_{i_1} x_{i_2} x_{i_3} x_{i_5} x_{i_6})(1+ x_{i_2} x_{i_3} x_
{i_4} x_{i_6} x_{i_7}).
\end{equation}  
Expanding out we see a term $x_{i_1} x^2_{i_2} x^2_{i_3} x_{i_5} x^2_{i_6} x_7$ corresponding to 
\begin{equation}
\chi(\mathcal{M}^{s, i_{6}}_Q(i_1 + 2i_2 + 2i_3 + i_5 +2i_6+ i_7)) = 1,
\end{equation} 
the unique contribution to $\chi(\mathcal{M}^{s, B}_{K(2)}(4, 6))$.
\end{exa}

Let us now move on to the general case of the base rings $\widetilde{R}_k$. We start with a result which says how to reconstruct the weight function $f_{\sigma_1 \cdots \sigma_l}$ from the tropical curve $h_{\sigma_1 \cdots \sigma_l}$. 
\begin{lem} The weight function attached to the tropical curve $h = h_{\sigma_1 \cdots \sigma_l}$ is given by
\begin{equation}\label{weightFunction}
f_{\sigma_1 \cdots \sigma_l} = \operatorname{Mult}_{Q}(h)\prod_{i, J}\left(\frac{(-1)^{\# J-1}}{\# J}(\# J-1)! \prod_{l \in J} u_{il}\right) x^{d_{\out}(h)},
\end{equation}
where the product is over all $i \in Q_0$ and $J \subset \{1, \dots, k\}$ for which $\mathfrak{d}_{i, J} \in \operatorname{Legs}(h)$, and 
\begin{equation}
d_{\operatorname{out}}(h) = \sum_{\mathfrak{d}_{i, J} \in \operatorname{Legs}(h)} (\# J)i.
\end{equation}
\end{lem}
\begin{proof} Suppose first that $l = 1$, so $h = h_{\sigma}$ for some $\sigma \in \mathfrak{S}^{\infty}$. If $\sigma$ first appears in $\mathfrak{S}^0_k$, then $h_{\sigma}$ is just one of the lines $\mathfrak{d}_{i, J}$, we have $\mult_Q(h_{\sigma}) = 1$ since there are no trivalent vertices at all, and $d_{\out} = (\# J)i$, $\ind(d_{\out}) = \# J$. Thus \eqref{weightFunction} holds by the definition of the operators $T_{i, J}$, \eqref{factoredOps}, and of $f_{\sigma}$, \ref{elemWeight}. Suppose now that $\sigma$ first appears in $\mathfrak{S}^i, k$ with $i \geq 1$, and $\operatorname{Parents}(\sigma) = \{\sigma_1, \sigma_2\}$. Then by induction and Lemma \ref{Commutator} we see that 
\begin{align}
\nonumber f_{\sigma} &= \mult_Q(h_{\sigma_1}) \mult_Q(h_{\sigma_2})\\
&\cdot \bra d_{1, \out}, d_{2, \out} \ket \prod_{i, J}\left(\frac{(-1)^{\# J-1}}{\# J}(\# J-1)! \prod_{l \in J} u_{il}\right) x^{d_{1,\out} + d_{2,\out}},
\end{align}
where the product is over all $i \in Q_0$ and $J \subset \{1, \dots, k\}$ for which $\mathfrak{d}_{i, J} \in \operatorname{Legs}(h_{\sigma_1}) \cup \operatorname{Legs}(h_{\sigma_2})$. Notice that $\operatorname{Legs}(h_{\sigma_1})\cap \operatorname{Legs}(h_{\sigma_2}) = \emptyset$ and $\operatorname{Legs}(h_{\sigma_1})\cup \operatorname{Legs}(h_{\sigma_2}) = \operatorname{Legs}(h_{\sigma})$. It is then clear by induction that 
\begin{equation}
d_{\operatorname{out}} = d_{1, \out} + d_{2, \out} = \sum_{\mathfrak{d}_{i, J} \in \operatorname{Legs}(h)} (\# J)i.
\end{equation}
On the other hand 
\begin{equation}
\mult_Q(h_{\sigma_1}) \mult_Q(h_{\sigma_2})\bra d_{1, \out}, d_{2, \out}\ket = \mult_Q(h_{\sigma}) 
\end{equation}
by the definition of multiplicity \eqref{QmultipCurve}.

For general $h_{\sigma_1\cdots\sigma_l}$ the same argument works, this time by induction on $l$. 
\end{proof}
\begin{lem} Let $\sigma_1 \prec \dots \prec \sigma_r$ be the maximal sequence of operators with slope $\mu$ in $\mathfrak{S}^{\infty}_k$. Then the product $\sigma_1 \circ \cdots \circ \sigma_r$ acts by
\begin{equation}
x_p \mapsto x_p \cdot f_p,
\end{equation}
where 
\begin{equation}\label{eff}
\log f_p = \sum_{h \in \mathcal{S}_{Q, k}(\mu)} \bra d_{\out}(h), i_p\ket f_{h}.
\end{equation}
\end{lem}
\begin{proof} As in the proof of Lemma \ref{Commutator} we write $\sigma_i$ as the exponential of some derivation in the noncommutaive Poisson algebra, $\exp(\{ \widetilde{\sigma}_i, \cdot\})$. Retaining only the first order corrections, the Baker-Camp\-bell-Hausdorff formula tells us that $\log \prod_i \exp(\{ \widetilde{\sigma}_i, \cdot\})$ is a sum of terms of the form 
\begin{equation}
\{\frac{1}{2}\{\dots \frac{1}{2}\{\frac{1}{2}\{\widetilde{\sigma}_{j_1}, \widetilde{\sigma}_{j_2}\},\widetilde{\sigma}_{j_3}\}, \dots, \widetilde{\sigma}_{j_l}\}, \cdot\}
\end{equation}
for $l \geq 1$ and $1 \leq j_1 < j_2 < \cdots < j_l$. However the formula \eqref{compose} in the proof of Lemma \ref{Commutator} tells us that this is an exact computation, i.e. only first order corrections actually happen. And again by \eqref{compose} and the definition of weight function \eqref{weightFunctDef}, we have
\begin{equation}
\frac{1}{2}\{\dots \frac{1}{2}\{\frac{1}{2}\{\widetilde{\sigma}_{j_1}, \widetilde{\sigma}_{j_2}\},\widetilde{\sigma}_{j_3}\}, \dots, \widetilde{\sigma}_{j_l}\} = f_{\sigma_{j_1}\cdots\sigma_{j_l}}.
\end{equation}
In other words there is a one to one correspondence between correction terms in the Baker-Campbell-Hausdorff formula and weight functions of (possibly disconnected) curves in $\mathcal{S}_{Q, k}(\mu)$.
\end{proof}
We will express $f$ in terms of our numbers $N^{\trop}_{Q, k}$. This computation should be compared with \cite{gps} Theorem 2.8 (although our case is simpler).
\begin{lem} 
The function $f_p$ in \eqref{eff} can be written as
\begin{equation}
f_p(t x) = \sum_{\emph{\w} : \mu(\emph{\w}) = \mu}\bra \emph{\w}, i_{p}\ket \frac{R_{\emph{\w}}}{|\Aut(\emph{\w})|} N^{\operatorname{trop}}_{Q, k}(\emph{\w})\,(t x)^{\emph{\w}}.
\end{equation}
\end{lem}
\begin{proof} For a fixed curve $h \in \mathcal{S}_{Q, k}(\mu)$ we find a weight vector $\w$ and sets $J_{qj}$ as in section \ref{curveConstruction}. We can rewrite \eqref{weightFunction} in terms of $\w$ as 
\begin{equation}
f_h = \operatorname{Mult}_Q(h)\prod^{s+S}_{q = 1} \prod^{l_q}_{j = 1}\left(\frac{(-1)^{w_{qj}-1}}{w_{qj}}(w_{qj}-1)!\prod_{r \in J_{qj}}u_{q r}\right)\,x^{\w}.
\end{equation}
Summing over all curves $h$ which give rise to the same weight vector $\textbf{w}$ and the same sets $J_{qj}$ we find a contribution to $\log f_p$ given by
\begin{equation}
\bra \w, i_{p} \ket N^{\textrm{trop}}_{Q,k}(\w)\prod^{s+S}_{q = 1} \prod^{l_q}_{j = 1}\left(\frac{(-1)^{w_{qj}-1}}{w_{qj}}(w_{qj}-1)!\prod_{r \in J_{qj}}u_{q r}\right)\,x^{\w}.
\end{equation}
Summing up over all $J_{qj}$ would then give 
\[
\bra \w,  i_{p}\ket N^{\trop}_{Q, k}(\w)\prod^{s+S}_{q = 1} \prod^{l_q}_{j = 1}\frac{(-1)^{w_{qj}-1}}{w^2_{qj}}\,t^{\w} x^{\w},
\]
but one can show that this overcounts curves by a factor $|\Aut(\w)|$.

\end{proof}
\begin{thm}\label{distanceOne} Let $i_p$ (for some $1 \leq p \leq s$) be a sink of $Q$ with precisely one source mapping to it, say $i_{\bar{p}}$ ($s + 1 \leq \bar{p} \leq s + S$). Then
\begin{equation}
\log \theta_{Q, \mu, \bar{p}}(t x) \equiv \sum_{\emph{\w} : \mu(\emph{\w}) = \mu} \bra i_p, \emph{\w}\ket \frac{R_{\emph{\w}}}{|\Aut(\emph{\w})|} N^{\operatorname{trop}}_{Q, k}(\emph{\w})\,(t x)^{\emph{\w}} \mod (t^{k+1}_1, \dots, t^{k+1}_{s+S}).
\end{equation} 
\end{thm}
\begin{proof} Reineke's theorem \ref{reinThm} gives
\begin{equation}
\theta_{\mu}(x_p) = x_p\cdot(\theta_{Q, \mu, \bar{p}}(t x))^{-1}.
\end{equation}
We also know 
\begin{equation}
\theta_{\mu}(x_p) \equiv x_p\,f_p\,(1 + \rho),\,\rho \in (t^{k+1}_1, \dots, t^{k+1}_{s+S}) 
\end{equation}
where $f_p$ is given by \ref{eff}. Therefore
\begin{align*}
\log \theta_{Q, \mu, \bar{p}}(t x) &\equiv -\log f_p \mod (t^{k+1}_1, \dots, t^{k+1}_{s+S})\\
&= \sum_{\w : \mu(\w) = \mu} \bra i_p, \w\ket \frac{R_{\w}}{|\Aut(\w)|} N^{\operatorname{trop}}_{Q, k}(\w)\,(t x)^{\w}.
\end{align*}
\end{proof}
We say that a vertex $i \in Q_0$ is a \emph{boundary vertex} ($i \in \partial Q$) if it is has valency $1$ in the undirected graph underlying $Q$. By making $Q$ larger if necessary, we can assume that the only boundary vertices of $Q$ are sinks. The result we just proved tells us how to compute $\theta_{Q, \mu, p}(t x)$ (to all orders) for all sources $p \in Q_0$ mapping to a boundary sink, in terms of certain tropical curves. Using the special feature that $Q$ sits in a tree (since the connected components of $\widetilde{K}(m)$ are infinite directed $m$-regular trees), we can propagate this calculation to an arbitrary sink of $Q$. 
\begin{thm} Let $Q \subset \widetilde{K}(m)$ be a finite connected subquiver, with only sinks as boundary vertices. Then for each source $i \in Q_0$, there exist distinct sinks $i^{pq}$ with $1 \leq p \leq P$, $1 \leq q \leq (m-1)^{p-1}$ such that
\begin{equation}
\log \theta_{Q, \mu, i}(x) \equiv \sum_{\emph{\w} : \mu(\emph{\w}) = \mu} \bra \varepsilon(i), \emph{\w}\ket \frac{R_{\emph{\w}}}{|\Aut(\emph{\w})|} N^{\operatorname{trop}}_{Q, k}(\emph{\w})\,(t x)^{\emph{\w}} \mod (t^{k+1}_1, \dots, t^{k+1}_{s+S}),
\end{equation}
with
\begin{equation}
\varepsilon(i) = i^{11} - (i^{21} + \dots i^{2(m-1)}) + \dots \pm (i^{P1} + \dots + i^{P(m-1)^{P-1}}).
\end{equation}
\end{thm}
\begin{proof}  Let us define recursively subsets $X_n, n \geq 1$ of the set of sources of $Q$, by
\begin{align*}
X_1 &= \{i \text{ is a source, } i \to j \in \partial Q\},\\
X_n &= \{i \text{ is a source, } i \to j \mid \text{ for all sources } i'\neq i \text{ with } i' \to j, i' \in X_{n-1}\}\cup X_{n-1}.
\end{align*}
Recall that we denoted by $S$ the number of sources. We claim that the set $X_S$ contains all the sources in $Q$. Arguing by contradiction, let $i_0$ be a source which is not contained in $X_S$. By the definition of $X_S$, we can pick any sink $j_0$ with $i \to j_0$, and find a source $i_1 \neq i$ with $i_1 \to j_0$ and $i_1 \notin X_{S-1}$. Notice that $i_1$ must map to a sink $j_1 \neq j_0$, otherwise $i_1$ would be a boundary vertex. Then in turn by the definition of $X_{S-1}$ we can find a source $i_2 \neq i_1$ with $i_2 \to j_1$ and $i_2 \notin X_{S-2}$. We must have $i_2 \neq i_0$ too otherwise $Q$ would contain an (unoriented) cycle. Proceeding by induction, given $i_n \notin X_{S-n}$ (with $n < S-1$) this maps to a sink $j_n \neq j_{n-1}$ (otherwise $i_n$ would lie on the boundary), and in fact $j_n$ is also distinct from all $j_0, j_1, \ldots j_{n-2}$ (otherwise $Q$ would contain an unoriented cycle); then by definition of $X_{S-n}$ we find $i_{n+1} \neq i_n$ with $i_{n + 1} \to j_n$ and $i_n \neq X_{S-n-1}$, $i_{n+1}$ also distinct from all $i_0, i_1, \ldots i_{n-1}$ (no cycles). We stop on reaching a source $i_{S-1}$. Thus we find a sequence of \emph{distinct} sources $i_0, i_1, \ldots i_{S-1}$, with $i_n \notin X_{S-n}$. In particular, this would say all the sources in $Q$ do not lie in $X_1$, which is clearly impossible. This proves the claim.
    
For a source $i$, we define $d(i, \partial Q)$ to be the least $n$ such that $i \in X_n$. We will prove our statement by induction on $d(i, \partial Q)$. If $d(i, \partial Q) = 1$ it reduces to Theorem \ref{distanceOne}. Suppose now the statement is known for all sources $j \in Q_0$ with $d(i, \partial Q) \leq P$ and choose $i \in Q_0$ with $d(i, \partial Q) = P + 1$. Then by the definition of distance $i$ maps to a sink $i^{11}$ such that all other sources $j_1, \dots, j_{m-1}$ mapping to $i^{11}$ (except $i$) satisfy $d(j_{q}, \partial Q) \leq P$. We apply Reineke's Theorem to $i^{11}$, finding
\begin{equation}
\theta_{Q, \mu}(x_{i^{11}}) = x_{i^{11}} (\theta_{Q, \mu, i}(tx))^{-1}\prod^{m-1}_{q = 1}(\theta_{Q, \mu, j_q}(tx))^{-1}.
\end{equation}
By our previous results, 
\begin{equation}
\theta_{Q, \mu}(x_{i^{11}}) = x_{i^{11}} f_{i^{11}},
\end{equation}
with 
\begin{equation}
f_{i^{11}} = \sum_{\w : \mu(\w) = \mu}\bra \w, i^{11}\ket \frac{R_{\w}}{|\Aut(\w)|} N^{\operatorname{trop}}_{Q, k}(\w)\,(t x)^{\w}.
\end{equation}
By induction, 
\begin{equation}
\log \theta_{Q, \mu, j_q}(tx) \equiv \sum_{\w : \mu(\w) = \mu} \bra \varepsilon(j_q), \w\ket \frac{R_{\w}}{|\Aut(\w)|} N^{\operatorname{trop}}_{Q, k}(\w)\,(t x)^{\w} \mod (t^{k+1}_1, \dots, t^{k+1}_{s+S})
\end{equation}
with 
\begin{equation}
\varepsilon(j_q) = j^{11}_q - (j^{21}_q + \dots + j^{2(m-1)}_q) + \dots \pm (j^{P1}_q + \dots + j^{P(m-1)^{P-1}}_q).
\end{equation}
Therefore
\begin{align}
\log \theta_{Q, \mu, i}(tx) &\equiv - \log f_{i^{11}} - \log \theta_{Q, \mu, j_1}(tx) - \dots - \log \theta_{Q, \mu, j_{m-1}}(tx)\\
&\equiv \sum_{\w : \mu(\w) = \mu} \bra \varepsilon(i), \w\ket \frac{R_{\w}}{|\Aut(\w)|} N^{\operatorname{trop}}_{Q, k}(\w)\,(t x)^{\w}
\end{align}  
with
\begin{align}
\nonumber &\varepsilon(i) = i^{11} - (j^{11}_1 + \dots +j^{11}_{m-1})\\
\nonumber &+ (j^{21}_1 + \dots + j^{21}_{m-1} + \dots + j^{22}_1 + \dots + j^{22}_{m-1} + \dots + j^{2(m-1)}_1 + \dots + j^{2(m-1)}_{m-1})\\
\nonumber & \vdots\\
& \pm (j^{P1}_1 + \dots + j^{P1}_{m-1} + \dots + j^{P2}_1 + \dots + j^{P2}_{m-1} + \dots + j^{P(m-1)^{P-1}}_1 + \dots + j^{P(m-1)^{P-1}}_{m-1})
\end{align}
where equivalence is modulo $(t^{k+1}_1, \dots, t^{k+1}_{s+S})$ as usual.
\end{proof}
Then Weist's Theorem yields the following.
\begin{cor}\label{contribGen} Let $\overline{d}$ be a fixed dimension vector for $K(m)$. Choose $Q \subset \widetilde{K}(m)$ large enough and $k \gg 1$ (depending on $\overline{d}$). Then
\begin{equation}
\chi(\mathcal{M}^{s, B}_{K(m)}(\overline{d})) = \sum_{[d]\sim \overline d} \sum^{s+S}_{p = s + 1}\exp\left( \sum_{\emph{\w} : \mu(\emph{\w}) = \mu} \bra \varepsilon(i_p), \emph{\w}\ket \frac{R_{\emph{\w}}}{|\Aut(\emph{\w})|} N^{\operatorname{trop}}_{Q, k}(\emph{\w})\,(t x)^{\emph{\w}}\right)[(tx)^d]
\end{equation}
where the first sum is over all equivalence classes $[d]$ of dimension vectors supported on $Q$ and compatible with $\overline{d}$.
\end{cor}
Finally, with these results in place, the GW/Kronecker correspondence gives back a comparison between the genuine tropical invariants $N^{\trop}(\w')$ of \cite{gps} and our ad hoc counts $N^{\trop}_Q(\w)$, at least in a special case. Indeed by the Kronecker/GW correspondence \ref{KronGW}, we have   
\begin{equation}
\chi(\mathcal{M}^{s, B}_{K(m)}(h a, h b)) = \exp\left(\frac{a}{m} \sum_{r \geq 1}\sum_{|P_a| = r a, |P_b| = r b} r N_{a, b}[(P_a, P_b)] x^{(ra, rb)}\right)[x^{(h a, h b)}]. 
\end{equation}
By Theorem \ref{tropicalThm}, the right hand side can be rewritten as
\begin{equation}
\exp\left(\frac{a}{m} \sum_{r \geq 1}\sum_{|P_a| = r a, |P_b| = r b} \sum_{\w'\sim (P_a, P_b)}\frac{R_{(P_a, P_b) | \w'}}{|\Aut(\w')|}N^{\trop}(\w') r x^{(ra, rb)}\right)[x^{(h a, h b)}].
\end{equation}
On the other hand for $Q$ and $k$ large enough the Euler characteristic can also be computed as
\begin{equation}
\sum_{[d]\sim (h a, h b)} \sum^{s+S}_{p = s + 1} \exp\left(\sum_{\w:\mu(\w) = \frac{a}{a+b}} \bra \varepsilon(i_p), \w\ket \frac{R_{\w}}{|\Aut(\w)|} N^{\operatorname{trop}}_{Q, k}(\w)(t x)^{\w}\right)[(tx)^d].
\end{equation}
Suppose now that $h = 1$. Since the vector $(a, b)$ is primitive, we must have the equality
\begin{align}
\nonumber \frac{a}{m}\sum_{|P_a| = a, |P_b| = b} \sum_{\w'\sim (P_a, P_b)} & \frac{R_{(P_a, P_b) | \w'}}{|\Aut(\w')|} N^{\trop}(\w') \\
&= \sum^{s+S}_{p = s + 1}\sum_{[\,\w\,] \sim (a, b)} \bra \varepsilon(i_p), \w\ket \frac{R_{\w}}{|\Aut(\w)|} N^{\operatorname{trop}}_{Q, k}(\w).
\end{align}
We do not know a direct proof of this equality, although we believe that one exists which exploits the comparison between the multiplicities $\mult_{Q, V} h$ and $\mult_V h$ as in \eqref{QmultipVert} and \eqref{multipVert2}.
\section{Connection with quiver gauge theories}\label{denef}
We will briefly explain the physical interpretation of $K(m)$ and how this picture (especially the paper of F. Denef \cite{denef}) gives a possible motivation for the constructions we have presented. Unfortunately the author is not an expert in the area, so our account will be very naive and imprecise. The reader should consult \cite{denef}.  

Let $S_1, S_2$ be two Lagrangian $3$-spheres in a compact Calabi-Yau threefold $X$, meeting transversely and positively in $m$ points, so for the intersection product (the DSZ product in this context) we have $\bra [S_1], [S_2]\ket = m$. In the terminology of \cite{denef} Section 3.1 $S_1, S_2$ are parton D3-branes. The generalised Kronecker quiver $K(m)$ with dimension vector $d = (d_1, d_2)$ arises in the study of the string theory on spacetime compactified on $X$ with $m$ open strings with boundaries on one of $d_1$ D-branes of type $[S_1]$ and one of $d_2$ D-branes of type $[S_2]$. 

The fundamental parameter in this theory is the string coupling constant $g_s$. For positive $g_s \approx 0$, and when the D-branes have small but nonvanishing phase difference and spacetime separation, the theory becomes a quiver quantum mechanics modelled on $K(m)$. In particular the Witten index of the theory can be computed as $\chi(\mathcal{M}_{K(m)}(d))$.  

A very different picture emerges for large coupling constant $g_s$. In this regime the BPS states for the theory become multi-centered, molecule-like configurations of $d_1$ ``monopoles" with charge $\textbf{Q}$ and $d_2$ ``electrons" with charge $\textbf{q}$, with DSZ product $\langle \textbf{Q}, \textbf{q} \rangle = m$ (i.e. the ``monopoles" have magnetic charge $m$, the ``electrons" have electric charge $1$). What Weist's Theorem \ref{weistThm} says in this regime is that we can compute the same Witten index by summing over all multi-centered BPS configurations with charges $\textbf{Q}_1, \dots, \textbf{Q}_{\ell_1}$ and $\textbf{q}_1, \dots, \textbf{q}_{\ell_2}$ such that the DSZ product $\langle \textbf{Q}_i, \textbf{q}_j\rangle$ is at most $1$ for $i = 1\, \dots, \ell_1, j = 1, \dots, \ell_2$ (i.e. such that each pair of interacting particles looks like a simple monopole-electron system, corresponding to $K(1)$). The $\widetilde{K}(m)$ constraint in this regime means that the splitting into charges $\textbf{Q}_i, \textbf{q}_j$ must be compatible with the original DSZ product $\langle \textbf{Q}, \textbf{q} \rangle = m$. 

For each of these multi-centered configurations, going back to $g_s \approx 0$ will give theories based on configurations of partons, with the same total Witten index. In other words one can compute the total Witten index by summing up over all the ways of splitting the boundary conditions for the open strings. Notice that so far we have ignored the framing, but this could easily be introduced by adding an additional parton D-brane $S$ to the discussion above. \\  
\begin{rmk}
The large $g_s$ viewpoint gives an interesting interpretation of Weist's gluing result \cite{weist} Corollary 5.28. Mathematically, in its simplest form, this says that if we have two representations $R', R''$ of $\widetilde{K}(m)$ with dimension vectors $d', d''$ and $\bra d', d'' \ket = 1$ we can glue them by identifying two sinks $j' \in R'_0, j'' \in R''_0$. The new dimension vector is $d = d' + d''$. Now from the $g_s \gg 0$ perspective we are simply superimposing our two special multi-centered configurations at two ``electrons". Weist's gluing corresponds to the statement that the total configuration we obtain is BPS, as long as the two multi-centered configurations behave mutually like a simple monopole-electron system, 
\begin{equation*}
\langle \sum_i \textbf{Q}'_i + \sum_j \textbf{q}'_j, \sum_i \textbf{Q}''_i + \sum_j \textbf{q}''_j\rangle = 1.
\end{equation*}
\end{rmk}
What we wish to retain from this physical picture is that passing to the universal cover for quiver representations corresponds to splitting either \emph{boundary conditions} (i.e. D-branes, partons) for small $g_s$ or particles (for large $g_s$) into a number of constituents. Then we can recover $\chi$ by summing up over all configurations of all possible types. The advantage of this physical point of view is that it suggests an analogy between  Weist's Theorem \eqref{weist} and Theorem \ref{tropicalThm}, i.e. in both cases we are computing our invariants (Witten indexes) by summing up over all boundary conditions (in other words it allows us to regard $[\,\tilde{d}\,]$ as specifying boundary conditions for open strings, while $\textbf{w}$ specifies ``boundary conditions", really legs, for tropical curves).

Finally we should mention that the special case $m = 2$ (with framing) has a physical interpretation as a certain $SU(2)$ Seiberg-Witten theory, as discussed in \cite{gmn} Section 2.2. The Kontsevich-Soibelman wall-crossing has been interpreted in the context of Seiberg-Witten theories by Gaiotto, Moore and Nietzke \cite{gmn}.  In the special case $m = 2$ the relevant identity is (using Kontsevich-Soibelman operators on $K(2)$) 
\[T_{1, 0} \circ T_{0, 1} = T_{0, 1} \circ T_{1, 2} \circ T_{2, 3} \cdots T^{-2}_{1,1}\cdots T_{3, 2} \circ T_{2, 1} \circ T_{1, 0}\]
which they interpret as going from strong coupling (the left hand side) to weak coupling (the right hand side). In Corollary \ref{contribSimple} we have seen which curves carry a contribution to the operator which represents one of the states of charge $(a, a+1)$.

\end{document}